\documentclass[preprint]{elsarticle}

\usepackage{amsmath}
\usepackage{ amssymb }
\DeclareMathOperator{\sech}{sech}
\usepackage{ dsfont }
\usepackage{mathtools}
\usepackage{color}
\usepackage{hyperref}
\usepackage{algorithm}
\usepackage{hhline}\usepackage[autostyle]{csquotes}
\definecolor{OrangeRed}{rgb}{0,0,0}
\definecolor{Magenta}{rgb}{0.9,0,0.9}
\definecolor{Black}{rgb}{0,0,0}
\newcommand{\magenta}[1]{{\color{Black}{#1}}}

\journal{Commun. Nonlinear Sci. Numer. Simulat.}

\begin{document}

\begin{frontmatter}

\title{Direct nonlinear Fourier transform algorithms for the computation of solitonic spectra in focusing nonlinear Schr\"{o}dinger equation}

\author[aston1,aston2]{A. Vasylchenkova$^*$}
\cortext[mycorrespondingauthor]{Corresponding author}
\ead[url]{vasylcha@aston.ac.uk}

\author[aston1]{J. E. Prilepsky}
\author[dm1,dm2]{D. Shepelsky}
\author[aston2,aston3]{A. Chattopadhyay}

\address[aston1]{Aston Institute of Photonic Technologies, Aston University, B4 7ET, Birmingham, UK}
\address[aston2]{System Analytics Research Institute, Aston University, B4 7ET, Birmingham, UK}
\address[dm1]{B. Verkin Institute for Low Temperature Physics and Engineering, Kharkiv 61103, Ukraine}
\address[dm2]{V. N. Karazin Kharkiv National University,  Kharkiv 61022, Ukraine}
\address[aston3]{Mathematics, Aston University, B4 7ET, Birmingham, UK}


\begin{abstract}
Starting from a comparison of some established numerical algorithms for the computation of the eigenvalues (discrete or solitonic spectrum) of the non-Hermitian version of the Zakharov-Shabat spectral problem, this article delivers new algorithms that combine the best features of the existing ones and thereby allays their relative weaknesses. Our algorithm is modeled within the remit of the so-called direct nonlinear Fourier transform (NFT) associated with the focusing nonlinear Schr\"{o}dinger equation. First, we present the data for the calibration of existing  methods comparing the relative errors associated with the computation of the continuous NF spectrum. Then each method is paired with different numerical algorithms for finding zeros of a complex-valued function to obtain the eigenvalues. Next we  describe a new class of methods based on the contour integrals evaluation for the efficient search of eigenvalues. After that we introduce a new hybrid method, one of our main results: the method combines the advances of contour integral approach and makes use of the iterative algorithms at its second stage for the refined eigenvalues search. The veracity of our new hybrid algorithm is established by estimating the convergence speed and accuracy across three independent test profiles. Along with the development of a
new approach for the computation of the eigenvalues, our study also
addresses the problem of computation of the so-called norming constants associated with the eigenvalues. We show that our formalism effectively amounts to accurate and fast enough computation of residues of the reflection coefficient in the upper complex half-plane of the spectral parameter.
\end{abstract}


\end{frontmatter}

\section{Introduction}\label{sec:intro}
Nonlinear Schr\"{o}dinger equation (NLSE) has traditionally played the de facto subservient modelling role over a wide range of topics that relate to the dynamical evolution of states and associated variables \cite{as,mg}. Such applications have been rampant in nonlinear physics, notably in photonics. In this article, our interest is in arriving at a generic master equation based formalism, structured around the NLSE, that governs the propagation of the slow-varying complex optical field envelope $q(z,t)$ along a single-mode lossless optical fibre \cite{mg,tpl17}:
\begin{equation}
i \frac{\partial q(t,z)}{\partial z} + \frac{1}{2}\,\frac{\partial^2 q(t,z)}{\partial t^2} + |q(t,z)|^2 q(t,z) = 0 .
\label{f:nse}
\end{equation}
For optical transmission-related problems, $z$ plays the role of the distance along the fibre while $t$ is the ubiquitous time variable, nomenclature that we will adhere to in the remainder of this article. It may be worth noting that similarly named quantities could have entirely different interpretations in other physical systems. In our study, we will only use suitably normalised dimensionless variables (both independent and dependent) to enable easy mapping of our results to other models within the same levels of description. The NLSE (\ref{f:nse}) has been explicitly written for the so-called {\it focusing} case, a term that is associated with \enquote{anomalous dispersion} in optical fibre studies \cite{mg}.

The celebrated work of Zakharov and Shabat \cite{zs72} revealed that Eq.~(\ref{f:nse}) belongs to a class of closed form integrable systems that can be completely solved by the inverse scattering transform method, subject to imposition of appropriate additional constraints on $q(t, z)$. In optical transmission literature, this method is popularly referred to as the {\it Nonlinear Fourier Transform} (NFT) method \cite{tpl17}, an allusion to the similarity with conventional Fourier transform that applies to linear PDEs. Such models have been shown to have robust applications in initial-value problems associated with nonlinear PDEs \cite{bc85}.
In our study, the NFT operation is landscaped to arrive at a complete set of nonlinear spectral data (\enquote{NF spectrum}) at a given spatial point $z=z_0$ by decomposing the known profile $q(t, z_0)$; here $q(t, z_0)$ effectively acts as the initial condition. The $z$-evolution of the individual NF spectral components is then decoupled and turns to be linear \cite{as,tpl17,zs72}. To find the space-time profile at a desired point $z=z_1$, we will need to solve the set of such linearised decoupled equations governing the NF spectrum evolution and then recover the spectrum distribution at $z_1$. $q(t, z_1)$ as a solution is uniquely evaluated by the reciprocal (inverse) NFT operation \cite{as,tpl17,zs72,bc85, yk14-1, stein}.

In this work, our primary aim is to compute the forward NFT assuming a bounded evolving initial condition, represented as follows:
\begin{equation} \int_{-\infty}^{\infty} \! |q(t, z_0)| \, dt < \infty.
\label{f:qcond1}
\end{equation}
This condition is automatically satisfied for each example profile studied in the remainder of this article. 

{\color{OrangeRed}{The spectral characteristics of the initial data consists of three parts: continuous data,
		defined in terms of the spectral (scattering) functions of 
		real spectral variable,
		eigenvalues (distinguished discrete values of spectral parameter), and respective norming constants associated with the eigenvalues.
		Accurate estimation of both parts of discrete NFT data is critically 
		important as they contain information about the bound states (solitons).
		Integrability of Eq.~(1) implies that if 
		$q(t,z)$ is a solution of Eq.~(1) with prescribed initial data
		$q(t,z_0)$ satisfying Eq.~(2), then $q(t,z)$ satisfies Eq.~(2)
		for all $z$ [8]. Moreover, the spectral data associated with $q(t,z)$ are independent of $z$, so hereafter, we remove the explicit reference of $q$ in $z$.}}


Inverse NFT is conventionally formulated over a subset of potentials $q(t)$, whose discrete spectrum satisfy the following additional constraints \cite{tf}:
\begin{enumerate}
	\item[(i)] complete lack of any real eigenvalue, and \label{item1}
	\item[(ii)] structurally simple eigenvalues. \label{item2}
\end{enumerate}
The second conditions could be technically regarded as a \enquote{soft} constraint but serves as a helpful guide during the numerical evaluation of these numbers. Also potentials
satisfying these assumptions are generic. They form an open dense set in the appropriate phase space, that relates to a subspace of (linear) operators originating from this functional space that could be directly applied to analyze stochastic communication problems, for example, in studying transmission over a quenched network.

We note here that all three aforementioned components of spectral data (continuous spectrum, eigenvalues, and norming constants) are interesting from the perspective of their usage in essentially nonlinearity-free optical transmission \cite{tpl17,yk14-1}. This is fundamentally important as nonlinearity is known to play the critical negative role of a \enquote{dampener} in high spectral efficiency optical communication systems, necessitating adequate measures to mitigate its impact on signal propagation, a key job for communication engineers \cite{ekw10,bmx16,cgk17}. This is where NFT-based methods could be highly beneficial, since the data transmission within the NFT-based framework occurs inside the NFT domain and thus is not hampered by nonlinearity-related pitfalls (recall that the NFT modes evolve linearly). While more complex technically, for distortionless communication chores, NFT-based signal processing could become a highly efficient alternative for the existing methods \cite{tpl17,yk14-1}.

A number of various communication systems based on the modulation of the different parts of NF spectrum have been proposed and studied recently. In the eigenvalue communications \cite{tm13} (name coined after the celebrated work by Hasegawa and Nyu \cite{hn93}), the complex soliton eigenvalues are used for modulation and transmission. The procedure is a natural generalization of the soliton-based methods \cite{mg}, the progenitor of NFT-based transmission. The further extension of this direction involves the utilization of norming constants together with eigenvalues \cite{yk14-1,hky14,dhg15,hyk16,hk16}. Our approach conforms closely with this methodology; the algorithms and signal processing NFT methods developed in our article target fast efficient computation of both the eigenvalues and norming constants. In addition to the eigenvalue based NFT transmission methods, we also refer to a group of methods that deal with the continuous part of the NF spectrum and non-solitonic NFT modes, either direct modulation based \cite{pdb14,lpp16} or within the so-called digital back-propagation framework \cite{tt13,wlp15}.
The computation of the continuous spectral functions is also addressed in our study. We use them to calibrate our methods.
Of late, there has been a significant progress in the experiments related to the NFT-based optical transmission \cite{dhg15,lpp16,b15} that could potentially relate to the theoretical framework that we propose to develop here. Finally, we note that Wahls and Poor recently proposed a faster algorithm for the computation of the NF spectrum \cite{wp13,wp15}, where they used a structure similar to the famous fast Fourier transformation (FFT) architecture \cite{ct65}. The computational time of this method grows as $n (\log n)^2$, where $n$ is the number of discretization points (samples). Our current study focuses exclusively on \enquote{conventional} methods, addressing the issue of their improvement, excluding for now the possibility of adapting our mechanism to undertake similar fast realizations as in \cite{ct65}.

The paper is organised as follows. Section \ref{sec:prem} introduces the mathematical formalism of the
NFT procedure and describes the quality metrics used in our further analysis.
Section \ref{sec:cont} introduces several NFT computation methods together with some of their modifications.
The accuracy of the spectral data are presented as follows:
\begin{enumerate}
	\item[(i)]  as a point-by-point variation of the spectral parameter $\xi$,
	\item[(ii)]  by using mean squared relative error,  and
	\item[(iii)] by using the energy associated with the continuous data of the reflection coefficient $r(\xi)$ (defined later).
\end{enumerate}
Section~\ref{sec:disc} contains the analysis of iterative and contour integration algorithms applied for the eigenvalue search, together with analysis of convergence and stability.
A new  hybrid approach for stable and accurate computation of the  eigenvalues  is detailed in this section.
Section~\ref{sec:norm} deals with the norming constants encompassing
two possible approaches leading to the calculation of the residue of the associated reflection coefficient
and several computational methods based on these approaches.
We show that our new approach is more accurate and stable (that is, converges faster) compared to the existing options. Conclusions and future directions are summarized in Section \ref{sec:concl}.

\section{Mathematical formulation, model signals, and performance metrics}\label{sec:prem}
\subsection{Formulation of the direct NFT operation and definition of NF spectrum quantities}\label{subsec:nft}
The NFT decomposition of a given pulse $q(t)$
is defined in terms of dedicated
solutions of
the Zakharov-Shabat system (ZSS) of ODEs \cite{zs72}
\begin{equation}
\frac{d}{dt}
\left(\begin{matrix}u_1(t, \xi)\\ u_2(t, \xi)\end{matrix}\right)
=\left(\begin{matrix}-i\xi&q(t)\\ -q^*(t)&i\xi\end{matrix}\right)
\left(\begin{matrix}u_1(t, \xi)\\ u_2(t, \xi)\end{matrix}\right),
\label{f:ZSode}
\end{equation}
where  $q(t)$, for  communication problems, is the signal to process (assuming the role of the ZSS potential decaying at the asymptotic limits $t \to \mp \infty$, in accordance with (\ref{f:qcond1})). $\xi$ is
the so-called spectral parameter which  can be understood as a nonlinear analogue of frequency. Asterisk in Eq. (\ref{f:ZSode}) and below denotes complex conjugates of corresponding quantities.
Apart  from the communication  and other NFT applications mentioned in Sec.~\ref{sec:intro}, the ZSS (\ref{f:ZSode}) is interesting by itself as it appears as a master equation in a coupled mode theory, describing, e.g. the scattering of waves in the fiber Bragg gratings \cite{k}, in electrical circuit related problems \cite{j82}, etc.

To retrieve the spectral data associated with the given profile $q(t)$, we fix the so-called Jost solutions 
$\Phi$ and $\Psi$
of Eq. (\ref{f:ZSode}) imposing
the asymptotic conditions at $t\to \pm \infty$:
\begin{equation}
\Phi(t,\xi)\equiv\left(\begin{matrix}\phi_1\\ \phi_2\end{matrix}\right)
\xrightarrow[t\rightarrow-\infty]{} \left(\begin{matrix}e^{-i\xi t}\\ 0\end{matrix}\right),
\quad
\Psi(t,\xi)\equiv\left(\begin{matrix}\psi_1\\ \psi_2\end{matrix}\right)
\xrightarrow[t\rightarrow \infty]{} \left(\begin{matrix} 0 \\ e^{i\xi t}\end{matrix}\right).
\label{f:asy}
\end{equation}
{\color{OrangeRed}{Using Eq.~(\ref{f:asy}), $\Phi(t,\xi)$ and $\Psi(t,\xi)$
		can be equivalently represented as solutions of Volterra integral equations
		(see, e.g., [8]), from which it is seen that they are 
		determined for all $\xi\in \overline{{\mathbb C}_+}$ 
		(the   upper complex half-plane closed by the real axis), are analytic 
		in $\mathbb C_+$ and continuous in $ \overline{{\mathbb C}_+}$.
		For all $\xi\in \overline{{\mathbb C}_-}$, we additionally define 
		$\hat \Psi(t,\xi)=\big(\psi_2^*(t,\xi^*), -\psi_1^*(t,\xi^*)\big)^T$.}}

The goal of the NFT pulse decomposition is to find the
continuous and discrete spectral quantities associated with $q(t)$.
{\color{OrangeRed}{Since $\Psi(t,\xi)$ and $\hat\Psi(t,\xi)$, considered 
		for $\xi\in \mathbb R$, constitute a fundamental system of solutions
		of (3), the representation of $\Phi(t,\xi)$ as a linear combination
		of these solutions
		introduces the  spectral (scattering) coefficients 
		$a(\xi)$ and $b(\xi)$ for $\xi\in \mathbb R$ by 
		\begin{equation}
		\Phi(t,\xi)=\hat\Psi(t,\xi)\:a(\xi) + \Psi(t,\xi)\:b(\xi), \qquad \xi\in {\mathbb R}
		\label{f:scattering}
		\end{equation}
		with 
		\begin{equation}\label{dets}
		a(\xi)=\det (\Phi(t,\xi),\Psi(t,\xi))
		\quad \text{and} \quad
		b(\xi)=\det (\hat\Psi(t,\xi),\Phi(t,\xi)).
	\end{equation}
		In view of (\ref{f:asy}), they can also be expressed by}}
\begin{equation}
a(\xi)=\lim_{t\rightarrow+\infty}\phi_1(t,\xi) e^{i\xi t},
\qquad b(\xi)=\lim_{t \rightarrow+\infty}\phi_2(t,\xi) e^{-i\xi t}.
\label{f:ab}
\end{equation}
They can also be characterized in terms of a single function, the reflection coefficient \cite{tf}:
\begin{equation}
r(\xi)=b(\xi)/a(\xi).
\label{f:r}
\end{equation}

{\color{OrangeRed}{The solitonic degrees of freedom are associated with the 
		discrete spectral data consisting of the set of complex-valued 
		eigenvalues $\{\xi_j\}$ of Eq. (3) that have positive
		imaginary parts (the set is finite due to condition (i) given 
		at the end of Sec.~1),
		together with complex-valued norming constants $\{c_j\}$.
		Since $\Phi(t\to -\infty, \xi)\to 0$ and 
		$\Psi(t\to \infty, \xi)\to 0$,
		for any $\xi$ with $\Im\xi>0$, it follows that an eigenvalue
		$\xi_j$ is characterized by the linear dependence of 
		$\Phi(t,\xi_j)$ and $\Psi(t,\xi_j)$, {\it i. e.} by the existence of 
		a non-zero constant $b_j\in \mathbb C$ such that 
\begin{equation}
\label{f:eigen}
\Phi(t,\xi_j) = \Psi(t,\xi_j)\: b_j.
\end{equation}
Hence the eigenfunction $\bigl(u^{(j)}_1(t),u^{(j)}_2(t)\bigr)^T$ of Eq.~(\ref{f:ZSode}),
associated with $\xi_j$, is given by 
$\bigl(u^{(j)}_1(t),u^{(j)}_2(t)\bigr)^T = \Phi(t,\xi_j) = \Psi(t,\xi_j)\: b_j$.

On the other hand, as shown in Eq.~(\ref{dets}), the eigenvalues can be 
equivalently characterized
as zeros of $a(\xi)$ in the upper half-plane: $a(\xi_j)=0$.
Generally, for potentials satisfying only Eq.~(\ref{f:qcond1}),
the constants ${b_j}$ are independent of the continuous scattering
functions, but if $\hat\Psi(t,\xi)$ 
(and thus $b(\xi)$) also admits analytic continuation into the upper half-plane 
(or at least in the domain $0<\text{Im} (\xi)<d$, for $d>0$, that is greater than the imaginary part of all eigenvalues), 
the parameters $\{b_j\}$ represent the 
values of $b(\xi)$ evaluated at $\xi=\xi_j$ (cf. (\ref{f:scattering})):}}
\begin{equation}\label{b-j}
b_j=b(\xi_j).
\end{equation}
In this case, the norming constants can be defined as
the residues of the reflection coefficient $r(\xi)$ at its poles $\{\xi_j\}$:
\begin{equation}
c_j = \mathrm{Res}[r(\xi)]|_{\xi = \xi_j} = \frac{b(\xi_j)}{a'(\xi_j)}
\label{f:res}
\end{equation}
{\color{OrangeRed}{(where we have assumed that the zeros of $a(\xi)$ are simple}}.
A sufficient condition ensuring analytic continuation of the above relates to a decay rate estimate for our profile $q(t)$:
\begin{equation}\nonumber|q(t)|<D e^{-d|t|}, \quad \text{for all}~t\in (-\infty,\infty),
\label{f:qcond2}
\end{equation}
for $D,d>0$. Particularly, for finitely supported $q(t)$ (which is the case
of the computational statement of the problem), $b(\xi)$ is analytic in the whole
plane and thus definition (\ref{b-j}) holds.

Individual soliton parameters can be directly extracted from the discrete  spectral data \cite{zs72}: the soliton amplitude is given by $2\text{Im} (\xi_j)$, and the soliton frequency  is $-2\text{Re} (\xi_j)$. The norming constant defines the remaining two soliton parameters: the center position of the individual soliton,
$$
j{\text{-th soliton centre position}}=-\frac{1}{2 \text{Im} (\xi_j)}\log \frac{|c_j|}{2 \text{Im} (\xi_j)},
$$and the solitonic phase that is proportional to the phase of the norming constant: $\varphi = -\text{arg} [i \, c_j]$. The aforementioned four real  parameters completely characterise each solitonic degree of freedom. Further details on the ZSS properties and soliton solutions can be found in Refs.~\cite{as,zs72,yk14-1} and \cite{tf}.

For our purpose, it will be useful to rewrite ZSS (\ref{f:ZSode}) for the wave envelope functions $\chi_{1,2}$ defined through the relations
\begin{equation}
\phi_1=\chi_1e^{-i\xi t}, \qquad
\phi_2=\chi_2e^{i\xi t}.
\label{f:phichi}
\end{equation}	
Then the ZSS for the envelope vector $X(t,\xi)=\big(\chi_1(t,\xi), \chi_2(t,\xi)\big)^T$ becomes
\begin{equation}
\frac{d}{dt}X(t,\xi)=\left(\begin{matrix}0&qe^{2i\xi t}\\ -q^*e^{-2i\xi t}&0\end{matrix}\right)X(t,\xi).
\label{f:ZSodesimplified}
\end{equation}
In terms of $X$, the spectral coefficients are given by
\begin{equation}a(\xi)=\lim_{t\to +\infty}\chi_1(t,\xi), \qquad b(\xi)=\lim_{t \to+\infty}\chi_2(t,\xi).
\label{f:absimpl}
\end{equation}
Since the initial conditions  for $X$ in Eq. (\ref{f:ZSodesimplified}) do not involve exponentials:
\[
(\chi_1(t,\xi), \chi_2(t,\xi))^T \to (1, 0)^T \qquad \text{ as} \ \ t \to -\infty,
\]
the definition of spectral coefficients via $X$
turns to be  convenient for some numerical methods described below.

\subsection{Model signals for test purposes}\label{subsec:signal}
Our algorithms will be tested against three independent test profiles. To ensure maximum possible variation, we have resorted to wide-spread model signals, where the analytical expressions for the spectral data can be written explicitly.

\begin{enumerate}
	\item[(i)] The over-soliton potential \cite{sy74} (or Satsuma-Yajima pulses) is given by
	\begin{equation}
	q_{\text{over}}(t)=A \sech t.
	\label{f:solpot}
	\end{equation}
	It is characterized by a single real amplitude parameter $A>0$.
	The associated spectral functions are as follows (all quantities with suffix \enquote{over} represent oversolitons):
	\begin{equation}
	a_{\text{over}}(\xi)=\frac{\Gamma^2 \left(\frac{1}{2}-i \xi \right)}{\Gamma \left(-A-i \xi +\frac{1}{2}\right) \Gamma \left(A-i \xi +\frac{1}{2}\right)},
	\label{f:sola}
	\end{equation}
	\begin{equation}
	b_{\text{over}}(\xi)=-\sin (\pi  A) \sech(\pi  \xi ),
	\label{f:solb}
	\end{equation}
	and
	\begin{equation}
	r_{\text{over}}(\xi)=-\frac{\sin (\pi  A) \, \text{sech}(\pi  \xi ) \, \Gamma \left(-A-i \xi +\frac{1}{2}\right) \, \Gamma \left(A-i \xi +\frac{1}{2}\right)}{\Gamma ^2\left(\frac{1}{2}-i \xi \right)},
	\label{f:solr}
	\end{equation}
	where $\Gamma(\ldots)$ is the Euler Gamma function.
	
	Depending on the value of  $A$, the discrete spectrum attributed to the oversoliton (\ref{f:solpot})
	consists of simple eigenvalues
	\begin{equation}
	\xi_k=(A-1/2-k)i,~k=0 \ldots \left[A-\frac{1}{2}\right],
	\label{f:soleigs}
	\end{equation}
	where $[\ldots]$ denotes the integer part. If $A$ is exactly half-integer, then $r(\xi)=0$ and the total energy is completely concentrated in
	the solitonic modes. The norming constant corresponding to the highest eigenvalue $\xi_0=(A-1/2)i$ is
	\begin{equation}
	c_{\text{over}}=i \Gamma (2 A)/ \Gamma^2 (A).
	\label{f:solres}
	\end{equation}
	More details on the  NFT properties of the profiles (\ref{f:solpot}) can be found in \cite{sy74}.
	
	\item[(ii)] The ZSS for the rectangular potential
	\begin{equation}
	q(t)=\begin{cases}
	A,~-L\leq t\leq L\\
	0,~\text{otherwise}
	\end{cases}
	\label{f:recpot}
	\end{equation}
	can also be solved analytically \cite{bo92,bct98}.
	The associated scattering coefficients are given by
	\begin{equation}
	a_{\text{rec}}(\xi)=e^{2i\xi L}\left(\cos \Big[ 2\sqrt{\xi^2+A^2} \, L \Big] -\frac{i \xi}{\sqrt{\xi^2+a^2}}\sin \Big[2\sqrt{\xi^2+A^2} \,L\Big]\right),
	\label{f:reca}
	\end{equation}
	\begin{equation}
	b_{\text{rec}}(\xi)=\frac{A}{\sqrt{\xi^2+A^2}}\sin \Big[2\sqrt{\xi^2+A^2} \, L \Big],
	\label{f:recb}
	\end{equation}
	and
	\begin{equation}
	r_{\text{rec}}(\xi)= \frac{A \exp(-2 i \xi L)}{i\xi -\sqrt{\xi^2+A^2}\cot \Big[ 2\sqrt{\xi^2+A^2} \, L\Big]}.
	\label{f:recr}
	\end{equation}
	The discrete eigenvalues $\{\xi_{rec}\}$ for the rectangle profile are given by the roots of the following transcendental equation
	\begin{equation}
	\tan \Big[ 2 \sqrt{A^2 + \xi_{rec}^2} \, L \Big] =\frac{\sqrt{A^2+\xi_{rec}^2}}{i \xi_{rec}}
	\label{f:recpole}
	\end{equation}
	for $\xi_{rec}$ in the upper half-plane. The norming constant for $\xi_{rec}$ is given by the expression
	\begin{equation}
	c_{\text{rec}}(\xi)=-\frac{i \left(A^2+\xi ^2\right) e^{-2 i L \xi }}{A \left(2 L \sqrt{A^2+\xi ^2} \cot \Big[ 2 \sqrt{A^2+\xi ^2} \, L \Big]-1\right)}.
	\label{f:recres}
	\end{equation}
	
	\item[(iii)] For our tests, we have also used the solitonic potential ($r(\xi)=0$, $\xi\in {\mathbb R}$) with unit amplitude and phase \cite{bct98}:
	\begin{equation}
	q_{\text{sol}}=\exp(-it)\sech(t).
	\label{f:shiftedpot}
	\end{equation}
	It has a single eigenvalue  $\xi_{\text{sol}}=0.5+0.5i$ with the associated norming constant
	$c_{\text{sol}}=i$.
	This potential allows us to check the behavior of our methods
	in the case of eigenvalues having a non-zero real part.
\end{enumerate}

\subsection{Accuracy and performance metrics}\label{subsec:accur}
Basically, any quantity that can be found analytically for ZSS solution can be used for the numerical NFT methods' accuracy assessment. For discrete $y$, we use relative error as an accuracy descriptor:
\begin{equation}
\epsilon=\frac{|y^{(\text{computed})}-y^{(\text{analytical})}|}{|y^{(\text{analytical})}|}.
\label{f:relerror}
\end{equation}
For a continuous spectral function $\varphi(\xi)$ ($\varphi$ can be $a$, $b$ or $r$),
we compare the analytical and computed values using the
mean squared relative error (MSRE):
\begin{equation}
\epsilon_{\varphi}=\frac{1}{N}\sum_{k=1}^N\frac{|\varphi_k^{\text{(computed)}}-\varphi_k^{\text{(analytical)}}|^2}{|\varphi_k^{\text{(analytical)}}|^2},
\label{f:integralerror}
\end{equation}
where $\varphi_k=\varphi(\xi_k)$.
For $\varphi^{\text{(analytical)}}(\xi)=0$,
we use the ordinary squared difference ${\big[|\varphi_k^{\text{(computed)}}-\varphi_k^{\text{(analytical)}}|\big]}^2$.

Following \cite{wp13,wp15}, we assess the performance of our algorithms using normalized runtime $\tau/n$, where $\tau$ measures the computer run time and $n$ is the number of discretization points (samples) defined over the finite processing interval  $t\in[-L \ldots L ]$.

For numerical evaluations, we used four nodes cluster with Intel(R) Xeon(R) CPU 2.50 GHz.

\section{Methods and algorithm calibration using continuous spectral functions}\label{sec:cont}
In this section, we present several NFT methods and compare their quality and performance
using the aforementioned model potentials over a range of test parameters. All of our methods use the potential truncation and its discretization over a finite time interval.

For the truncation of model potentials and discretization, $q(t)$ is analytically represented to define samples within the interval $t\in[-L, L]$.
The interval is divided into $n$ equal subintervals of length $\Delta t=2L/n$, where the ${m}$-th subinterval is $t\in[t_m-\Delta t/2, t_m+\Delta t/2]$. Outside of the interval $[-L,L]$, the potential is assumed to be exactly zero. For all methods considered below, our signal (the ZSS potential) is approximated by a constant value along a single step: $q_m=q(t_m)$.

The vector of wave functions $\Phi(t)=(\phi_1(t), \phi_2(t))^T$ is fixed by imposing the initial conditions
on the left edge of the truncation interval ($t=-L$): according to (\ref{f:asy}), $\Phi(-L)=(e^{i\xi L}, 0)^T$.
The vector of  envelopes  $X(t)=(\chi_1(t), \chi_2(t))^T$ is fixed correspondingly as $X(-L)=(1,0)^T$.
The spectral functions $a(\xi)$ and $b(\xi)$ are defined on the right edge of the processing interval ($t=L$):
according to (\ref{f:ab}) and (\ref{f:absimpl}), $a(\xi)=\phi_1(L)e^{i \xi L}$ and $b(\xi)=\phi_2(L)e^{-i \xi L}$, or,
alternatively,
$a(\xi)=\chi_1(L)$ and $b(\xi)=\chi_2(L)$.

\subsection{Transfer matrix methods for NF spectrum computation}\label{subsec:meth}
To propagate the incident wave towards the end of the processing interval, many approaches use the
transfer matrix $T_m$ for propagating the ZSS solution over a single ${m}$-th discretization step, inside which the potential is considered as a constant, $q_m=\text{const}$:
\begin{equation}
\Phi_{m+1}=T_m\Phi_m.
\label{f:evolution}
\end{equation}
Performing the consequent iterations of Eq. (\ref{f:evolution}) from $m=1$ to $m=n$, we eventually find the desired values of the Jost functions at the end of the processing interval and compute the NFT parameters $a(\xi)$ and $b(\xi)$.

\begin{enumerate}
	\item[(i)] For the Bofotta-Osborn (BO) method \cite{tpl17,bo92,bct98}, we define $T_m$ {\color{OrangeRed}{evaluating the matrix exponential of the matrix in ZSS with constant potential $q_m$:}}
	\begin{equation}
	T_m^{\text{(BO)}}=\exp\left[\Delta t\left(	\begin{matrix} -i\xi & q_m \\ -q_m^* &  i\xi \end{matrix}	\right)\right] \! ,
	\label{f:BOT1}
	\end{equation}
	which can be evaluated explicitly:
	\begin{equation}
	T_m^{\text{(BO)}}=\left(\begin{matrix}	\cosh \kappa \Delta t - i\xi / \kappa \sinh \kappa \Delta t & q_m/\kappa \sinh \kappa \Delta t\\
	-q_m^*/\kappa \sinh \kappa \Delta t & 	\cosh \kappa \Delta t + i\xi / \kappa \sinh \kappa \Delta t\end{matrix}	\right) ,
	\label{f:BOT}
	\end{equation}
	with $\kappa=\sqrt{-|q_m|^2-\xi^2}$.
	
	\item[(ii)] For the Ablowitz-Ladik (AL) method, we use the normalized discretization of ZSS \cite{wp13,wp15,yk14-2} {\color{OrangeRed}{and apply Euler method, then substituting $1\pm i\xi \Delta t $ with $\exp(\pm i \xi \Delta t)$, we get the transfer matrix:}}
	\begin{equation}
	T_m^{\text{(AL)}}=\frac{1}{\sqrt{1+\Delta t^2 |q_m|^2}}\left(
	\begin{matrix}	
	e^{-i\xi \Delta t} & \Delta t q_m\\
	-\Delta t q_m^* & 	e^{i\xi \Delta t}
	\end{matrix}	\right).
	\label{f:ALT}
	\end{equation}
	The AL method with the norming factor $1/\sqrt{1+\Delta t^2 |q_m|^2}$ provides  higher stability and accuracy compared to that without the normalisation \cite{wp13}.
\end{enumerate}

In this article, we also introduce and study two novel modified versions of BO and AL algorithms,
for which the similar approaches are applied to the ZSS for the envelope functions (\ref{f:ZSodesimplified}). Evolution over each step $\Delta t$ is again performed using the transfer matrices:
\begin{equation}
X_{m+1}=T_mX_m.
\end{equation}
\begin{enumerate}
	\item[(iii)] {\color{OrangeRed}{For the modified BO method, applying the matrix exponential to the matrix of ZSS for the envelope functions defined in Eq.~(13), one can get the transfer matrix in the form:}}
	\begin{equation}
	T_m^{\text{(BOmod)}}=\left(\begin{matrix}\cos|q_m\Delta t| &\sin|q_m\Delta t|e^{i(\theta_{q_m}+2\xi t)}\\ -\sin|q_m\Delta t|e^{-i(\theta_{q_m}+2\xi t)}&\cos|q_m\Delta t|\end{matrix}\right),
	\label{f:BOmodT}
	\end{equation}
	where $\theta_{q_m}$ is $\textrm{arg}[q_m]$.
	
	\item[(iv)] In the case of the modified AL method, we have
	\begin{equation}
	T_m^{\text{(ALmod)}}=\frac{1}{\sqrt{1+\Delta t^2 |q_m|^2}}\left(
	\begin{matrix}	
	1 & q_m\Delta t e^{2i\xi t}\\
	-q_m^*\Delta t e^{-2i \xi t} & 1
	\end{matrix}\right).
	\label{f:ALmodT}
	\end{equation}

	\item[(v)] Finally, for the (non-modified) Crank-Nicolson (CN) method \cite{yk14-2,wp13}, the transfer matrix entering Eq.~(\ref{f:evolution}) is given by
	\begin{equation}
	T_m^{\text{(CN)}}=(I-\frac{\Delta t}{2}P_{m+1})^{-1}(I+\frac{\Delta t}{2}P_m),
	\label{f:CNT}
	\end{equation}
	where $I$ is the $2\times2$  identity matrix and
	\begin{equation}
	P_m=\left(
	\begin{matrix}	
	-i\xi &  q_m\\
	- q_m^* & 	i\xi
	\end{matrix}
	\right).
	\label{f:CNT1}
	\end{equation}
	
\end{enumerate}

Notice that in the NFT related works (see e.g.~\cite{yk14-2,bct98}) the algorithms for the solutions of ordinary ODEs (the Runge-Kutta scheme) were also studied in the application of the ZSS analysis. The Runge-Kutta fourth-order algorithm (RK) for the solution of ZSS is better applicable to the envelope system Eq. (\ref{f:ZSodesimplified}), as in this case the rapid oscillations of the Jost functions (for the region of $\xi$ with a large real part) are included into the effective potential functions. We do not describe the RK scheme here as it is quite standard (see e.g. \cite{bct98}) but below we present a comparative analysis of the RK algorithm with the algorithms mentioned above.

\subsection{Results for the continuous spectral functions}\label{subsec:cont}
\begin{figure}
	\centering
	\includegraphics[width=0.4\textwidth]{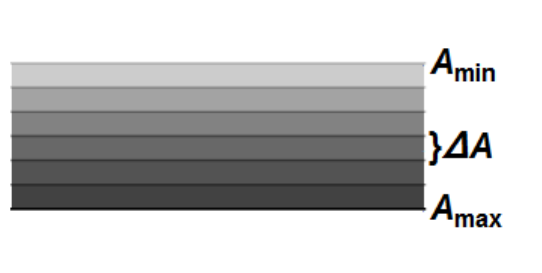}
	\caption{Transparency scale for amplitude variation in grey tones; progressively darker shades indicate higher amplitudes}
	\label{fig:shades}
\end{figure}

Now we perform a comparative accuracy analysis of the described methods using the three descriptors: $a(\xi)$, $b(\xi)$, and $r(\xi)$, Eqs. (\ref{f:ab}),  (\ref{f:r}). For the real values of $\xi$ the MSRE (\ref{f:integralerror}) is used for the accuracy assessment.
We also address the behavior of the NFT methods in analyzing the dependence of the method's accuracy on the variation of amplitude: in all the following figures, the amplitude changes are depicted using the transparency scale of the corresponding curves, see Fig.~\ref{fig:shades}, i.e. the curves for different amplitudes are plotted by the colored areas changing from an almost transparent part (the lowest amplitudes) to an almost solid line (the highest amplitude). Captions to the plots provide information about the range and step of amplitude variation.

\begin{figure}
	\centering
	
	\begin{tabular}{cc}
		\includegraphics[width=0.45\textwidth]{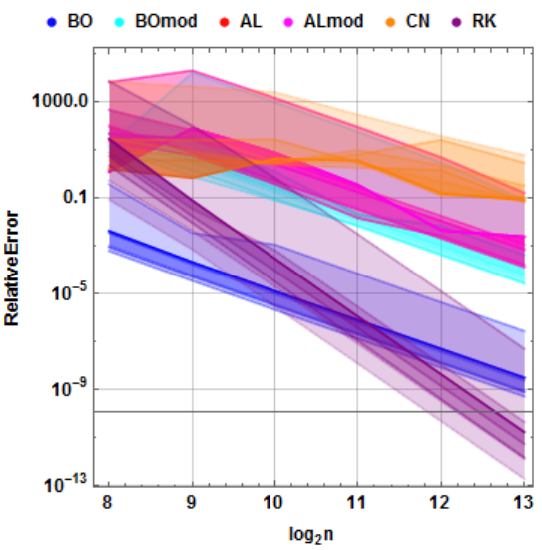}&
		\includegraphics[width=0.45\textwidth]{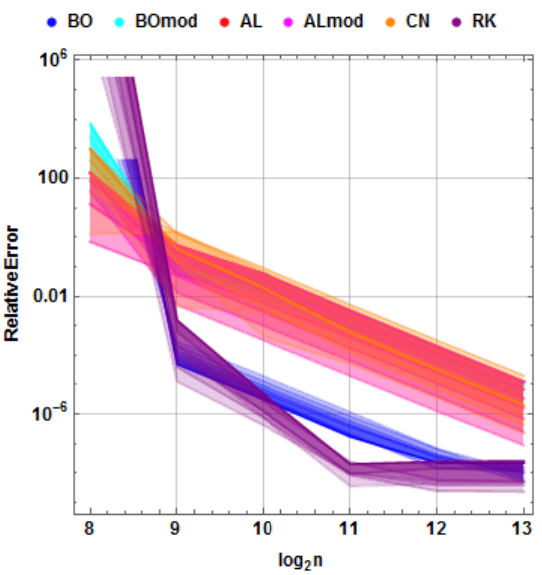}\\
		{\scriptsize{a) $a(\xi)$ error for $q_{\text{over}}$}}&
		{\scriptsize{b) $b(\xi)$ error for $q_{\text{over}}$}}\\
		\includegraphics[width=0.45\textwidth]{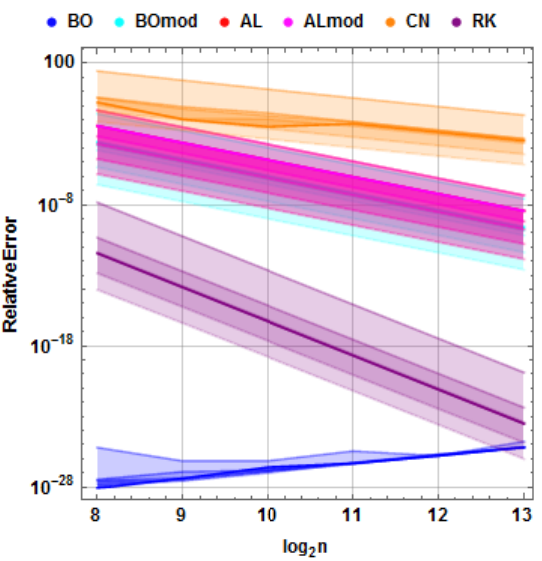}&
		\includegraphics[width=0.45\textwidth]{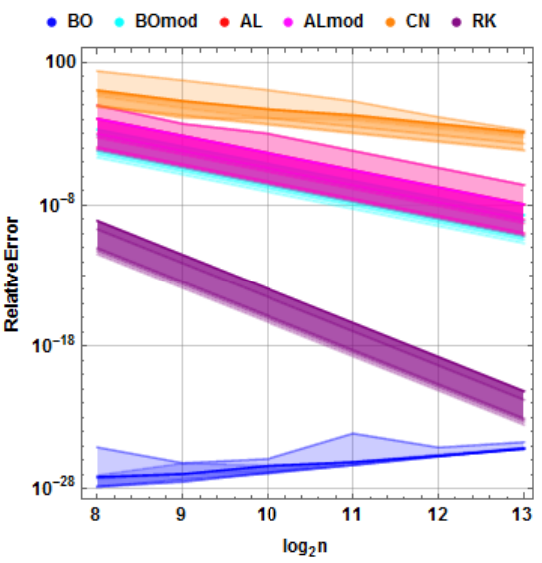}\\
		{\scriptsize{c) $a(\xi)$ error for $q_{\text{rec}}$}}&
		{\scriptsize{d) $b(\xi)$ error for $q_{\text{rec}}$}}\\
	\end{tabular}
	\caption{MSRE for the computation of the NF spectrum parameters (the specific spectral parameters are given in the captions) as a function of the number of discretization points $n$, evaluated for different NFT methods, and shown in a)--b) panels: for the over-soliton potential $q_{\text{over}}$, where amplitude changes in the range $[2.25, \, \ldots, \, 5.25]$ with the increment step $\Delta A=0.5$ and $L=30$; c)--d) panels: for the rectangular potential $q_{\text{rec}}$, where amplitude changes in the range $[2, \, \ldots, \, 5]$  with the increment step $\Delta A=0.5$ and $L=1$.}
	\label{fig:n}
\end{figure}

Our analysis confirms that the BO method gives the best accuracy among all methods studied, see Fig.~\ref{fig:n}.  The AL and the modified BO algorithms display similar behavior with the change of amplitude $A$ and of the number of points $n$, whereas the CN algorithm came up with the worst accuracy and convergence rate. Both AL and BO methods have the same convergence (inclination of the curves on logarithmic scale plots), implying that they all have the same order of accuracy. This conclusion complies with the results summarized in \cite{tpl17} and earlier studies. As expected, the fourth order RK method converges faster. This method can be better than BO
for big enough values of $n$. At the same time, for a smaller number of discretization points, the RK method's MSRE is excessively large. The BO method shows the weakest dependence on the amplitude variation, whereas the RK is the most sensitive to it (the error increment can reach several orders of magnitude in the range of amplitudes that we used for our plots). We also note that the BO method in application to the rectangular potential  gives the solution, which coincides with explicit analytical expression, so that the main source of errors here is the computational error in evaluation of $\cosh$ and $\sinh$ from (\ref{f:BOT}). This offsets the surprising increase in error with increase in the number of points for this particular potential and method.

The energy embedded into the continuous NF spectrum can also be a convenient quantity for methods' accuracy assessment: in particular, it was used as a qualifying metric in~\cite{bct98}. This energy is defined through the spectral functions $a(\xi)$ or $r(\xi)$ as follows
\begin{equation}
\varepsilon=-\frac{1}{\pi}\int_{-\infty}^{\infty}\log |a(\xi) |^2 d\xi=\frac{1}{\pi}\int_{-\infty}^{\infty} \log (1+|r(\xi)|^2).
\label{f:energy}
\end{equation}
(Recall that $|a(\xi)|^2+|b(\xi)|^2=1$ for $\xi\in {\mathbb R}$).

\begin{figure}
	\centering
	\begin{tabular}{cc}
		\includegraphics[width=.62\textwidth]{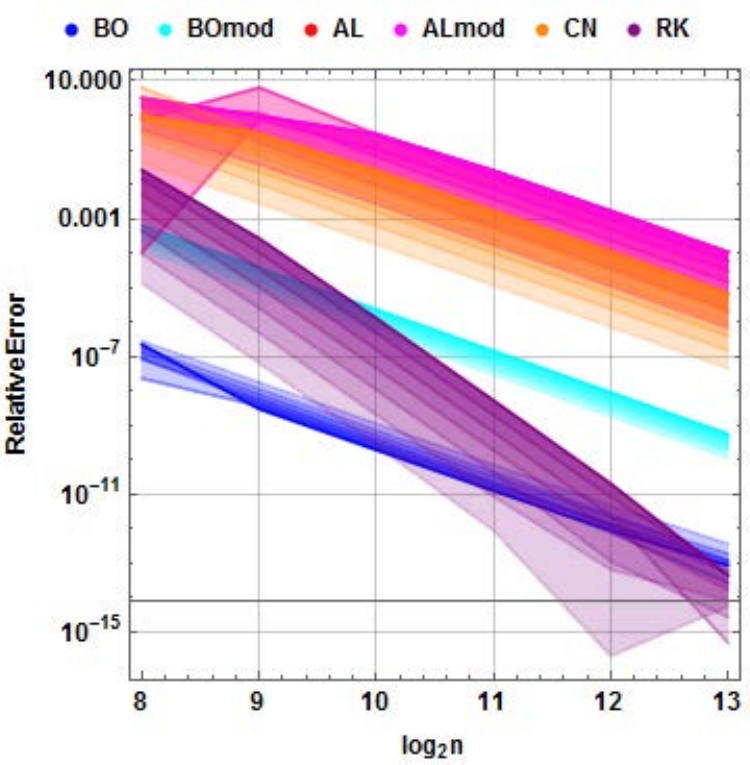}&
		\begin{tabular}[b]{c}
			\includegraphics[width=.37\textwidth]{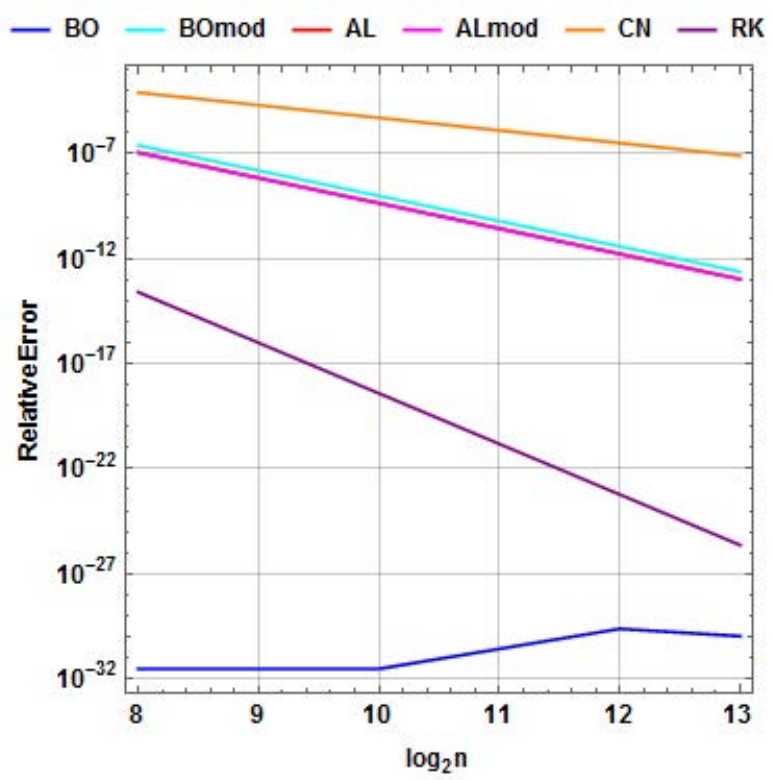}\\
			{\scriptsize{b) $q_{\text{rec}}$, $L=1$, $A=\pi/2$}}\\
			\includegraphics[width=.37\textwidth]{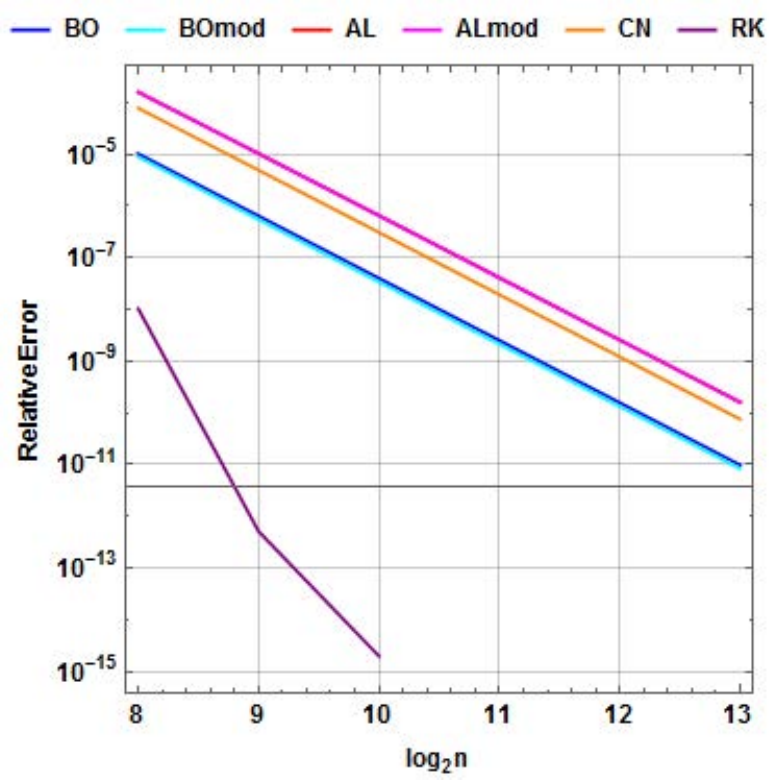}\\
			
		\end{tabular}\\
		{\scriptsize{a) $q_{\text{over}}$, $L=30$}}&	{\scriptsize{c) $q_{\text{sol}}$, $L=30$}}\\
	\end{tabular}
	\caption{Relative error $\epsilon$ (\ref{f:relerror}) in the continuous spectrum energy $\varepsilon$ (\ref{f:energy}) versus the number of discretization points $n$ for the different test potentials. The variation with the potential amplitude change is shown for over-soliton signal only (left pane), where $A$ changes in the range $[2.25, \, \ldots, \, 5.25]$ with the increment step $\Delta A=0.5$. }
	\label{fig:enerrorn}
\end{figure}

We analyze the relative error in the calculation of energy (\ref{f:relerror}) versus the number of points and amplitude, see Fig.~\ref{fig:enerrorn}.
The energy analysis shows qualitatively similar result as the MSRE analysis does, which confirms its correctness for the NFT methods' accuracy assessment. The RK algorithm converges more rapidly than all other methods, especially for the pure solitonic potential, where the double precision numbers are not enough to find out the difference between the analytical and computed energy values. At the same time, this test reveals yet another disadvantage of the RK method: it is extremely slow in comparison with all other methods (see Table \ref{t:cont}). The second worst in terms of time consumption is the CN method, whereas both the AL and the modified BO are similar in terms of computational time consumption and are the fastest among the all discussed methods. These methods are approximately two times faster than the ordinary BO algorithm, but the latter has an impressively high accuracy.
The results for the amplitude dependence when using the energy as a metric, Fig.~\ref{fig:enerrorn}, are similar to our findings when the MSRE for the  NFT continuous data was used, Fig.~\ref{fig:n}.

{\color{OrangeRed}{\begin{table}
	\caption {Runtimes (in seconds) of continuous energy evaluation for different NFT algorithms for $n=2^{13}$ for different test profiles}\label{t:cont}
	\centering
	\def\arraystretch{1.2}
	\begin{tabular}{|c|c|c|c|c|c|}
		\hline
		\multicolumn{6}{|c|}{NFT method}\\
		\hline
		BO& BOmod &AL &ALmod&CN&RK\\
		\hline
		\multicolumn{6}{|c|}{$q_{\text{over}}$, $L=30$, $A=5.25$}\\
		\hline
		125.28&71.8&70.11&71.33&232.31&1552.99\\
		\hline
		\multicolumn{6}{|c|}{$q_{\text{rec}}$, $L=1$, $A=\pi/2$}\\
		\hline
		177.43&104.21&121.03&124.99&307.55&2799.38\\
		\hline
		\multicolumn{6}{|c|}{$q_{\text{sol}}$, $L=30$}\\
		\hline
		89.99&49.52&56.14&58.32&145.67&1308.59\\
		\hline
	\end{tabular}
\end{table}
}}

\begin{figure}
	\centering
	\begin{tabular}{cc}
		\includegraphics[width=0.45\textwidth]{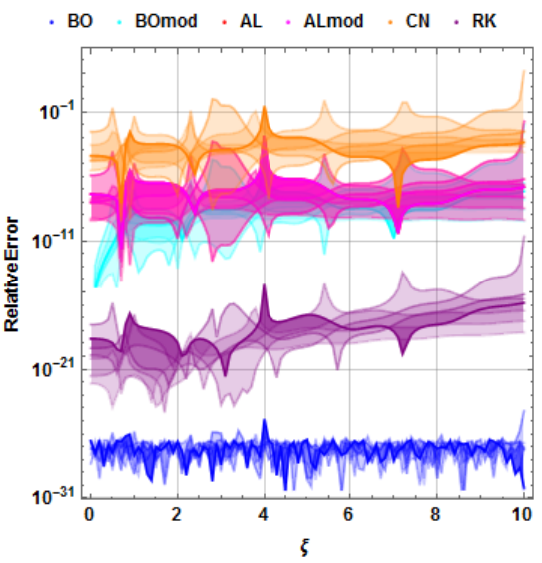}&
		\includegraphics[width=0.45\textwidth]{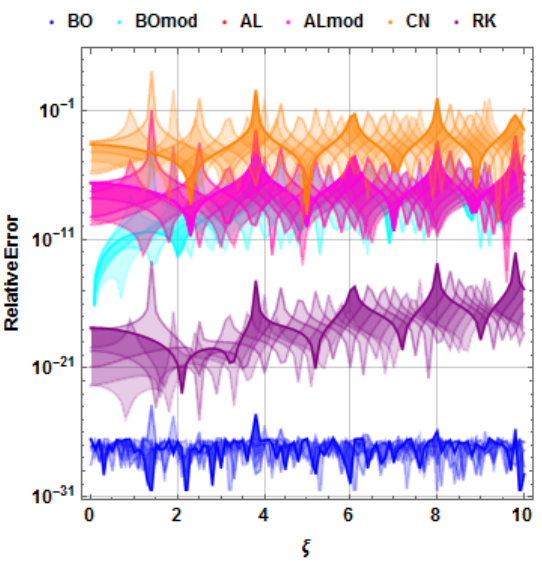}\\
		{\scriptsize{a) $a(\xi)$ for $q_{\text{rec}}$}}&
		{\scriptsize{b) $b(\xi)$ for $q_{\text{rec}}$}}\\
		\includegraphics[width=0.45\textwidth]{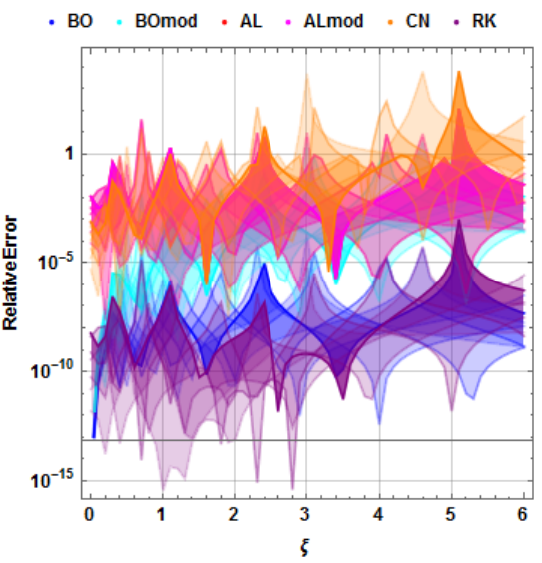}&
		\includegraphics[width=0.45\textwidth]{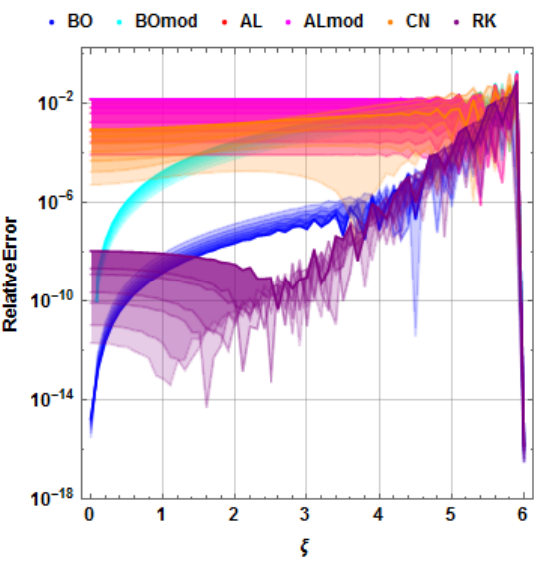}\\
		{\scriptsize{d) $a(\xi)$ for $q_{\text{over}}$}}&
		{\scriptsize{e) $b(\xi)$ for $q_{\text{over}}$}}\\
	\end{tabular}
	\caption{Relative error versus $\xi$ {\color{OrangeRed}{(with step size $\Delta\xi=0.1$)}} and amplitude, panels  a)--b): for $q_{\text{rec}}$ with $L=1$ and amplitude changes in the range $[2, \, \ldots, \, 5]$ with the increment step $\Delta A=0.5$; panels c)--d): for $q_{\text{over}}$ with $L=20$ and amplitude changes in the range $[2.25, \, \ldots, \, 5.25]$ with the increment step $\Delta A=0.5$. {Profiles are sampled with $n=2^{10}$}. }
	\label{fig:errorxi}
\end{figure}

In order to analyze the NFT methods stability in dependence on the nonlinear frequency bandwidth, we investigate
the accuracy of our methods along the nonlinear frequency $\xi$ axis point-by-point, see Fig.~\ref{fig:errorxi}.
{\color{OrangeRed}{As concluded before, BO and RK methods typically show higher accuracy. In estimating the accuracy dependence of the frequency, we observe that the numerical error deviates dramatically for the above-mentioned methods, while for AL, ALmod, BOmod and CN methods the numerical error is a lot more stable. This effect is more pronounced for the over-soliton potential.}}
This tendency does not change significantly when we are tuning the amplitude of the potential. The CN method demonstrates the worth accuracy among the all methods studied.

Now let us compare our findings with the previous results. The accuracy assessment of computing the continuous spectral data was done in \cite{bct98}.  The BO and RK methods were compared there with regard to the continuous spectrum energy computation: the convergence of the methods was studied and their runtimes were analyzed. It was shown that for smooth solitonic potentials, the RK method was better than the BO method, but the authors \cite{bct98} attributed this finding to the properties of the CPU architecture used for their computations. For the rectangular potential, which has sharp edges, the RK method was shown to perform slower than the BO one. 
	{\color{OrangeRed}{We have generally observed that in terms of
		runtime and resulting accuracy, the BO method typically outperforms the other
		options; we have also noticed the excessively large runtime required by the RK
		methods in comparison to any transfer-matrix based approach of the type shown in Eq.~(\ref{f:evolution}).}}. 
	The authors of \cite{wp13} and \cite{yk14-2} presented the detailed description of various NFT algorithms, but they were mainly focused on the discrete eigenvalues computation accuracy. The authors of \cite{wp13} also proposed the fast implementation of the AL algorithms based on the FFT-type ideology for the matrix product computation and multipoint fast polynomial evaluation, and this study was continued in \cite{wp15} for the periodic NFT potential. The comparison of the \enquote{conventional} non-fast NFT methods accuracy for the case of periodic potentials and discrete eigenvalues were also presented in \cite{kpl16}. Typically, the qualitative behavior of the accuracy as a function of $n$  followed the scenario described above in this section.

\section{Computation of  eigenvalues}\label{sec:disc}
The solitonic eigenvalues of ZSS under the assumptions listed above are defined as zeros of the analytical extension of $a(\xi)$, see (\ref{f:ab}), into the upper half-plane $\mathds{C}^+$ of the complex $\xi$  plane. In this section we present two principally different approaches to the computation of location of solitonic eigenvalues and then compare their accuracy, stability and performance.
\subsection{Iterative methods for eigenvalue search}\label{subsec:iter}
For the  eigenvalues computation, the most popular option among the NFT related works is to apply some iterative scheme to identify the complex zero(s) of a function (namely, $a(\xi)$ from Eq. (\ref{f:ab})) in the case considered \cite{yk14-2}. In  \cite{wp15,yk14-2,bct98}, the dependence of the NFT method performance on a particular iterative scheme usage was somewhat overlooked and an arbitrarily chosen method was usually employed without a particular motivation for the choice. In our paper, we fill this gap and compare the existent iterative schemes: we analyze the convergence and rapidity of the computational methods in application to the spectral data computation using the test profiles from Subs.~\ref{subsec:signal}.

Traditionally, the most common iterative approaches are the secant and Newton-Raphson (NR) methods. Assuming an initial guess for the location of the zero $x_0$ of some function, say $f(x)$, the consecutive iteration scheme for the NR method is given by
\begin{equation}
x_{k+1}^{\text{(NR)}}=x_k-\frac{f_k}{f'_k},
\label{f:NR}
\end{equation}
where  $f_k:=f(x_k)$. This method has a quadratic convergence rate \cite{j03} ({\color{OrangeRed}{numerical method is said to have convergence rate $p$, if $|x_{k+1}-x|\leq C |x_k-x|^p$,}} here and below the orders of convergence are given under the assumption that all roots are simple). The main disadvantage of the NR method is the necessity to know the value of the function derivative at each iteration step.
For the purpose of brevity, we introduce a shorthand notation for the so-called divided differences:
\begin{equation}\label{f:divdif}
f[x_1, x_2 ] =\frac{f(x_2)-f(x_1)}{x_2-x_1},
\end{equation}
which can be {\color{OrangeRed}recursively} generalized for an arbitrary number of arguments:
\magenta{
	\begin{equation}
f[x_k, \ldots, x_{k+m}]=\frac{f[x_{k+1},\ldots, x_{k+m}]-f[x_{k}, \ldots, x_{k+m-1}]}{x_{k+m}-x_k}.
\label{f:divdif1}
\end{equation}
}
In the secant method, the expression for the derivative in each iteration is swapped over to the divided difference, leading to:
\begin{equation}
x_{k+1}^{\text{(secant)}}=x_k-\frac{f_k}{f[x_{k-1}, x_k]}=x_k-f_k\frac{x_k-x_{k-1}}{f_k-f_{k-1}}.
\label{f:secant}
\end{equation}
The convergence rate of the secant method is approximately 1.618, so it is worse than that for the NR method, but a single step computation using the secant method can be faster since it does not require  computing the derivatives.
Sidi \cite{s08} generalised the idea of the derivative approximation: the function derivative $f'_k$
can be replaced by  the derivative of a fitting polynomial $p(x)$ of degree $j$:
\magenta{
\begin{equation}
p'_j(x)=f[x_{k-1},x_{k}]+\sum_{i=2}^{j}f[x_{k-i}, ...x_{k}]\prod_{l=1}^{i-1}(x_k-x_{k-l}).
\label{f:Sidi1}
\end{equation}
}
The next iteration is given by
\begin{equation}
x_{k+1}^{\text{(Sidi)}}=x_k-\frac{f_k}{p'_j}.
\label{f:Sidi}
\end{equation}
In our study, we use a cubic polynomial approximation in (\ref{f:Sidi1}), i.e. $j=3$. For this particular case, the convergence rate of the method is $\approx 1.93$ \cite{s08}.

Steffensen's method \cite{js68} uses the following iterative formula:
\begin{equation}
x_{k+1}^{\text{(Steffensen)}}=x_k-\frac{f_k^2}{f(x_k+f_k)-f_k}.
\label{f:St}
\end{equation}
It allows us to reach the convergence rate $2$, same as that for the NR method.

The Muller method \cite{m56} has an advantage in that it allows us to find complex roots from a real initial guess. Defining the auxiliary quantities
\begin{equation}
\begin{aligned}
&w=f[x_{k-1}, x_k]+f[x_{k-2}, x_k]-f[x_{k-2}, x_{k-1}], \qquad g=f[x_{k-2}, x_{k-1}, x_k],\\
&d=\text{max}\left[w-\sqrt{w^2-4f_k g}, ~w+\sqrt{w^2-4f_k g}\right],
\end{aligned}
\end{equation}
(the maximum is determined by comparing the absolute values), the iteration step of the Muller method is given by
\begin{equation}
x_{k+1}^{\text{(Muller)}}=x_k-2\cdot f_k/d.
\label{f:Muller}
\end{equation}
The order of convergence for this method is approximately $1.84$ \cite{j03}, which is better than that for the secant method.

\begin{figure}
	\centering
	\includegraphics[width=0.4\textwidth]{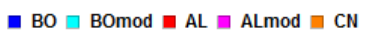}
	\begin{tabular}{ccc}
		\includegraphics[width=0.32\textwidth]{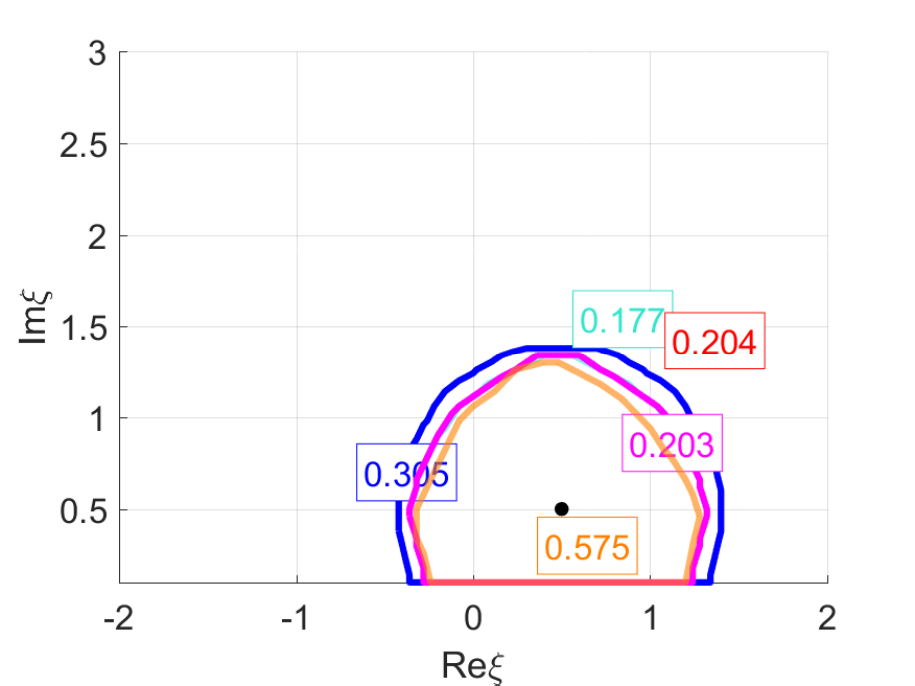} &
		\includegraphics[width=0.32\textwidth]{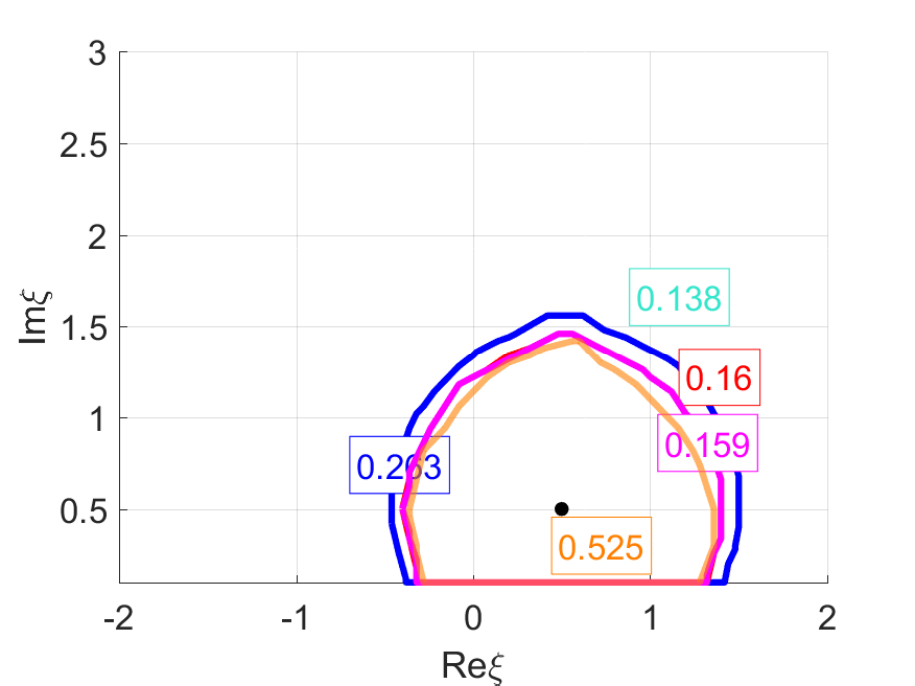}&
		\includegraphics[width=0.32\textwidth]{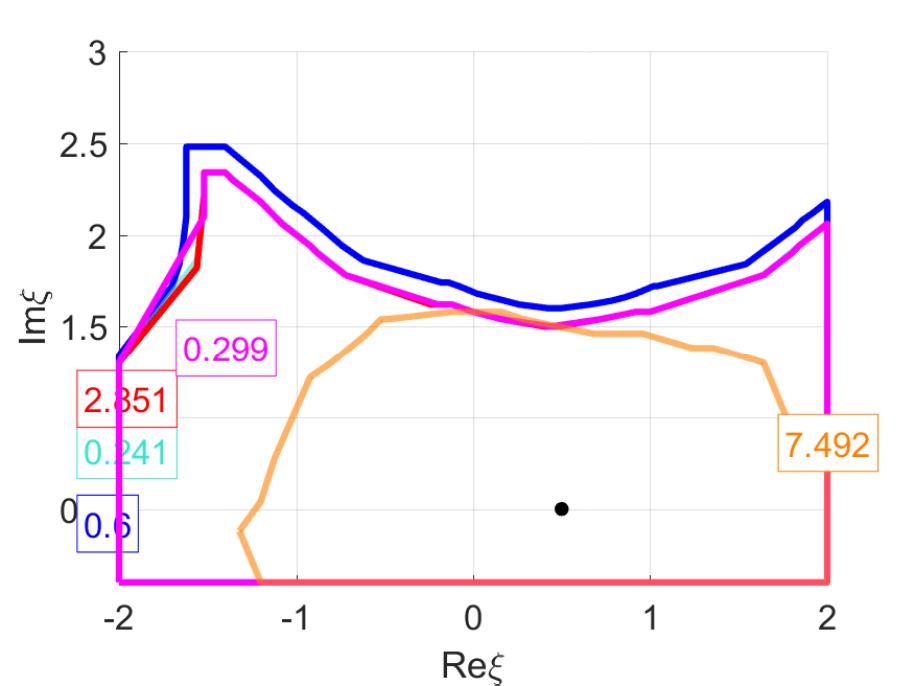} \\
		{\scriptsize{a) Newton-Raphson}}&{\scriptsize{b) secant}}&{\scriptsize{c) Sidi}}\\
	\end{tabular}
	\begin{tabular}{cc}
		\includegraphics[width=0.32\textwidth]{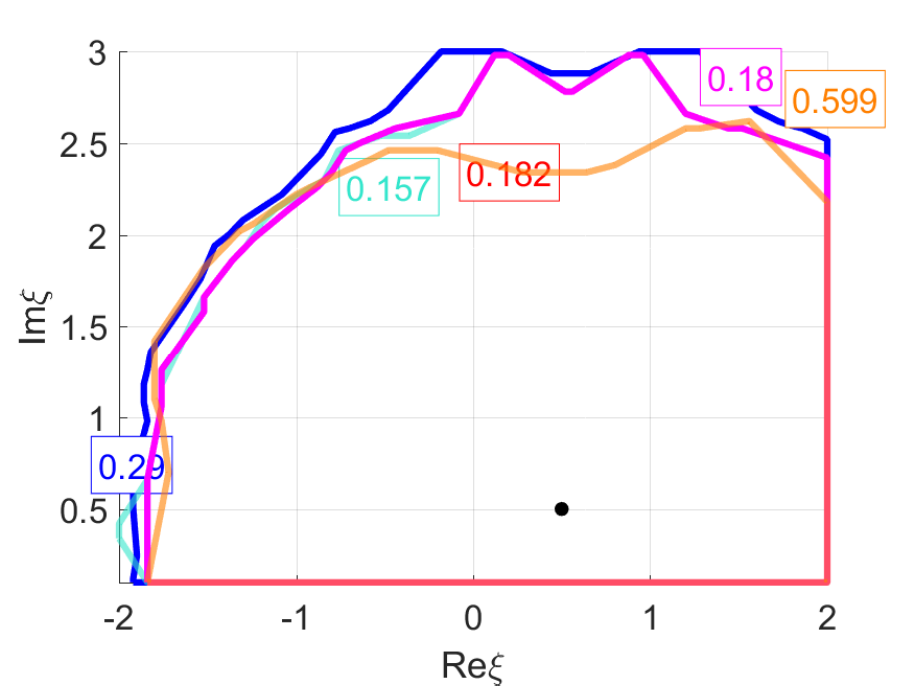}&
		\includegraphics[width=0.32\textwidth]{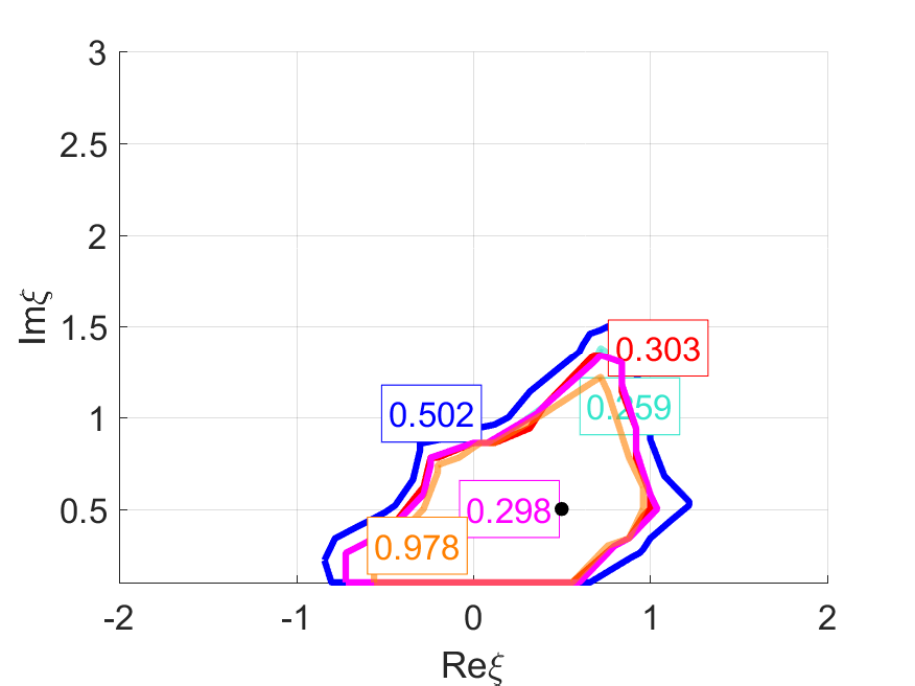}  \\
		{\scriptsize{d) Muller}}&{\scriptsize{f) Steffensen}}\\
	\end{tabular}
	\caption{The borders of the regions of initial assumptions for the zero approximation, from which the corresponding iterative algorithm reaches its correct value (marked as black point) with less than 1\% relative error.	
		Computation was performed for the soliton with phase factor potential $q_{\text{sol}}$ given by Eq. (\ref{f:shiftedpot}) for the values $n=2^{10},~L=30$. The digits in the legends near each curve identify the average runtime for the computation of an eigenvalue, when the initial guess point was taken inside the respective regions. Perfect vertical and horizontal edges of some basins mean that regions of convergence extend out of ranges where we make these probes.}
	\label{fig:poolsol}
\end{figure}

\begin{figure}
	\centering
	\includegraphics[width=0.4\textwidth]{Figlegend}
	\begin{tabular}{ccc}
		\includegraphics[width=0.32\textwidth]{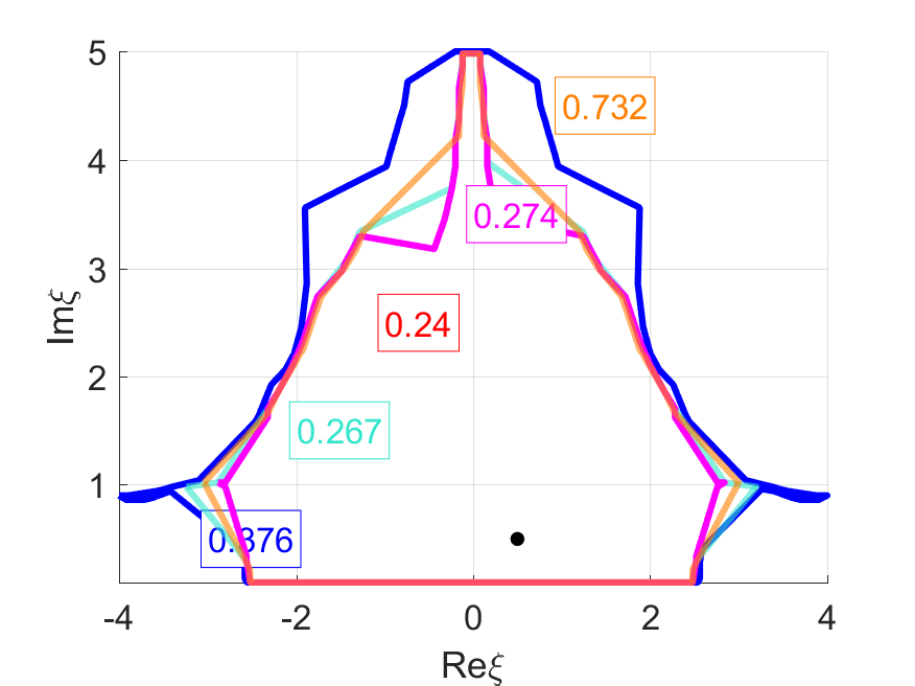} &
		\includegraphics[width=0.32\textwidth]{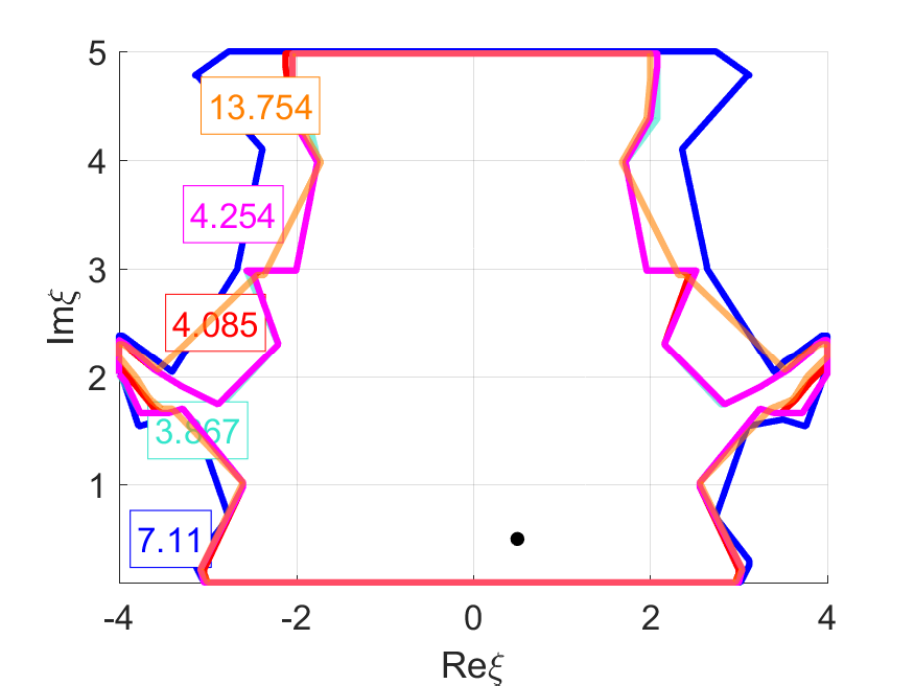}&
		\includegraphics[width=0.32\textwidth]{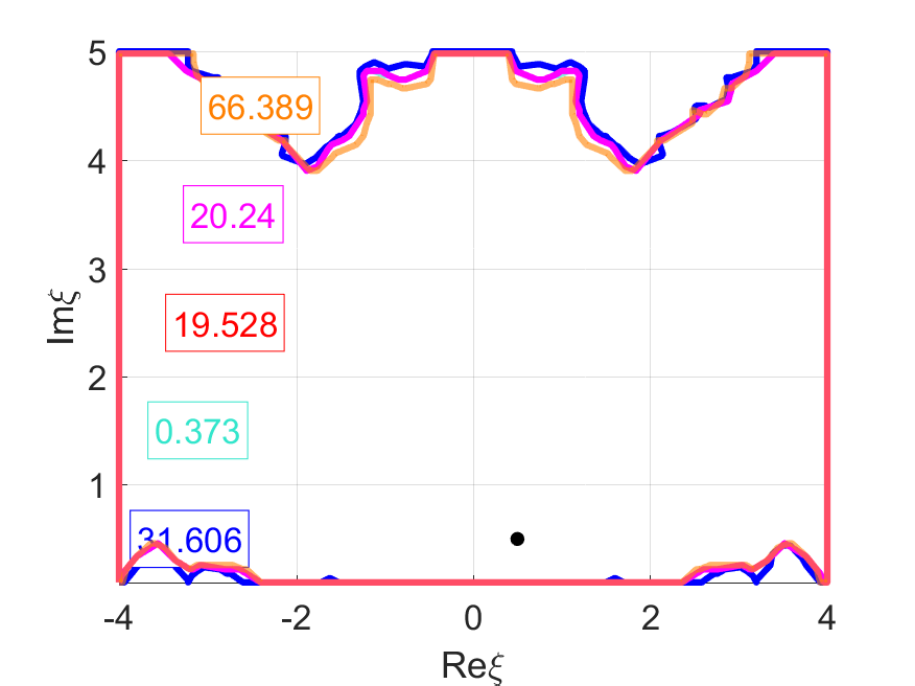} \\
		{\scriptsize{a) Newton-Raphson}}&{\scriptsize{b) secant}}&{\scriptsize{c) Sidi}}\\
	\end{tabular}
	\begin{tabular}{cc}
		\includegraphics[width=0.32\textwidth]{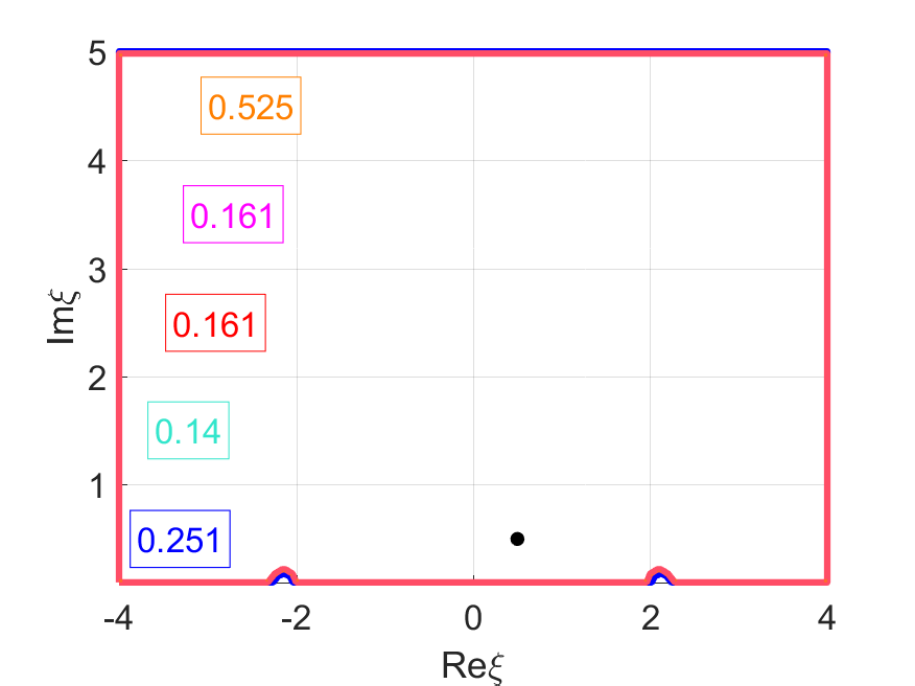}&
		\includegraphics[width=0.32\textwidth]{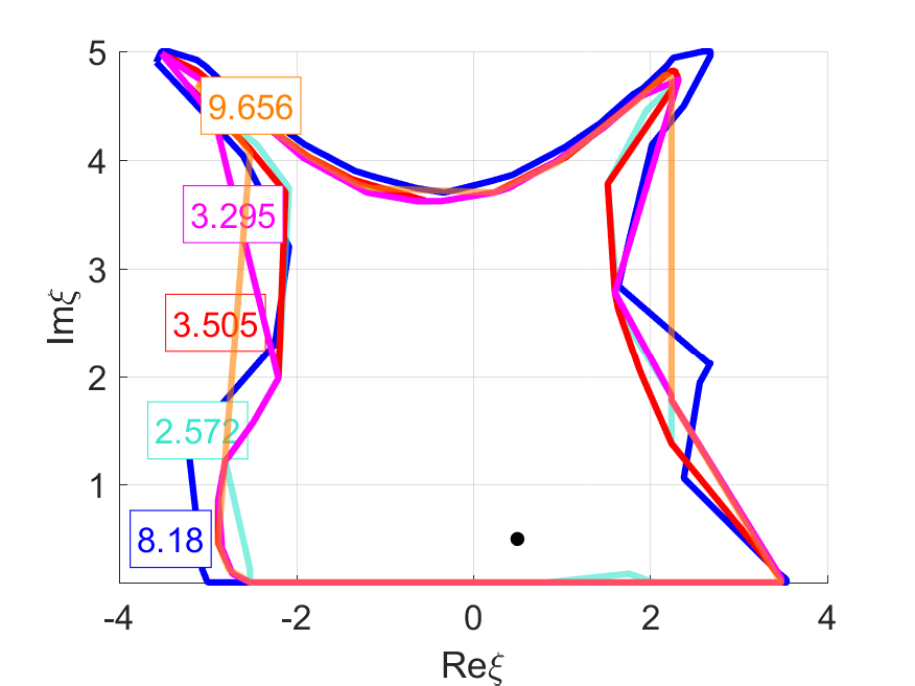}  \\
		{\scriptsize{d) Muller}}&{\scriptsize{f) Steffensen}}\\
	\end{tabular}
	\caption{The borders of the regions of initial guesses for zero approximation, from which the corresponding iterative algorithm reaches its best value (marked as black point) with less than 1\% relative error. Computation was performed for rectangular potential $q_{\text{rec}}$ ($n=2^{10},~L=1,~A=\pi/2$). The digits in the legends near each particular line identify the average runtime for the computation of an eigenvalue, when the initial guess point was taken inside the respective regions. Perfect vertical and horisontal edges of some basins mean that region of convergence extends out of ranges where we make these probes. }
	\label{fig:poolrec}
\end{figure}

All iterative algorithms applied for the eigenvalues computation require good initial guess. In order to understand how the choice of the initial  value influences the result of the  eigenvalues search, we investigate
the convergence of all iterative schemes in dependence on the initial guess point value: Fig.~\ref{fig:poolsol} contains the results referring to the solitonic potential with phase $q_{\text{sol}}$, and in Fig.~\ref{fig:poolrec} we depict the results for the rectangular potential $q_{\text{rec}}$ with $A=\pi/2$ where only a single eigenvalue is present. On these plots we show the border of the regions in the complex plane of spectral parameter $\xi$ (marked with the closed lines of different colour), where the relative error of zero location (estimated by using (\ref{f:relerror})) is less than $0.01$. 	{\color{OrangeRed}{We run the iterative algorithm until it reaches the pre-set precision in the difference between the function values for consequent iterations or until it exceeds the pre-set number of iterations (these pre-sets were correspondingly $10^{-10}$ and $10^3$).}} We also add to each line on the plots the corresponding average runtime that the computation of the eigenvalue takes when the initial guess point is positioned inside the respective regions.

The result of our analysis of the aforementioned five iterative methods combined with different NFT algorithms can be summarized as follows: the NR and secant methods, used  in \cite{wp15,yk14-2}, have smaller regions of convergence than those of the Sidi and Muller methods (the Muller method shows the largest convergence basin for all example profiles used), whereas the runtimes for all four approaches are similar. In the case of the Sidi algorithm, we observe a better convergence, but with significantly higher runtime, especially when the iterative method is coupled with the AL and CN NFT algorithms. The Steffensen iterative method shows even smaller region of convergence, which is worse than that for the NR and secant methods.

The drawback of iterative methods in the application to real transmission problems (when the position of solitonic eigenvalues is usually not known apriori) is that we cannot predict the computational runtime since we cannot estimate the number of iterations required to find the zero point with a satisfactory accuracy. Moreover, the methods can fail to converge at all, so that some additional precautions have to be taken. The previous works devoted to the eigenvalue search algorithms assumed iterative schemes for estimating discrete eigenvalues. The main results referred to particular features of the NFT method but not necessarily to a particular iterative scheme. Authors of \cite{yk14-2} studied the one-soliton and multisoliton cases for the AL, Euler and RK methods (also using the CN method for some cases, but evidently revealing the weakness of the CN method). They also found out that all studied root-search techniques resulted in similar accuracy regimes for the NFT data computed.
According to \cite{bct98}, the RK algorithm can converge faster than the BO, but the authors there used the grid search for the location of eigenvalues, and this resulted in a relatively high overall runtime. It was found  that in the case of rectangular potential, the RK method failed in the correct computation for the total number of zeros (the authors used the total phase increment along the $\text{Re} (\xi)$ axis for this purpose).
Our analysis of the iterative algorithms reveals that they are not sufficiently stable and manageable for the  eigenvalue computations in realistic applications. This fact motivates us to seek for principally different methods and options applicable for the location of eigenvalues.

What we have not studied here, though, are the class of methods that involve matrix diagonalization, e.g. the Fourier collocation method \cite{yk14-2}. The idea of the Fourier collocation is the decomposition of the ZSS problem in the Fourier series and the subsequent diagonalization of a specific complex-valued non-Hermitian block matrix, recasting the zero search as an eigenvalue problem. According to \cite{yk14-2}, this method can provide a good accuracy (its accuracy is spectral, in contrast to the methods studied in our current paper), but it has so far been adapted to find the eigenvalues only. Usually it also requires a considerable runtime ($\sim n^3$) to reach the result.

\subsection{Contour integration methods for eigenvalues search}\label{subsec:contour}
In this subsection, we present a new alternative technique for eigenvalue computation. Delver and Lyness in \cite{dl67} presented an approach for the location of the complex zeros of nonlinear functions based on the contour integral evaluation, see also \cite{ksb99} for more resent results on this approach. This method will be identified as DL (or ordinary DL) in the rest of this text. Within the method of \cite{dl67}, we start from the relation:
\begin{equation}
\frac{1}{2\pi i}\int_C z^p \frac{f'(z)}{f(z)}dz=\sum z_i^p,
\label{f:DL1}
\end{equation}
that allows us to make up the root-search scheme that would be applicable for finding all zeros $z_i$ of a function $f(z)$ inside the closed contour $C$ drawn in the complex plane. Within this approach, the zeros $z_i$ will emerge as the roots of some specially constructed polynomial. First, setting $p=0$ in (\ref{f:DL1}) gives us the total number of zeros, $N$, located inside the contour $C$. Next, evaluating (\ref{f:DL1}) for $p=1 \ldots N $ we can readily find the sums of $z_i^p$ up to $p=N$:
\begin{equation}
s_p=\sum z_i^p, \qquad p=1\ldots N.
\label{f:DL2}
\end{equation}
Having found the set of $s_p$, we can write down the equation system to evaluate the so-called Newton's identities, $\sigma_p$:
{\color{OrangeRed}{
		\begin{equation}
		\begin{aligned}
		\sigma_1&=-\sum z_i:~s_1+\sigma_1=0\\
		\sigma_2&=z_1z_2+z_2z_3+...+z_{n-1}z_N:~s_2+s_1\sigma_1+2\sigma_2=0\\
		...&\\
		\sigma_N&=(-1)^Nz_1z_2...z_N:~s_N+s_{N-1}\sigma_1+...+s_1\sigma_{N-1}+N\sigma_N=0.
		\end{aligned}
		\label{f:DL3}
		\end{equation}
}}
This system of equations can be solved recurrently using the values of the set \{$\sigma_i$\} obtained in the previous iteration round:
\begin{equation}
\label{f:DL4}
\sigma_p=\frac{1}{p}\left(\sum_{j=1}^{p-1}s_j\sigma_{p-j}+s_p\right).
\end{equation}
Next, using the Newton's identities we can construct the following polynomial:
\begin{equation}
P(z)=z^N+\sigma_1z^{N-1}+\sigma_2z^{N-2}+\ldots+\sigma_{N-1}z+\sigma_N.
\end{equation}
The polynomial $P(z)$ has exactly the same roots as the initial function $f(z)$. Therefore, using any polynomial root-finding technique, the desired set of roots of function $f(z)$, $z_i$, can be estimated. If the initial function has multiple roots, then they will be presented repeatedly along with the set of polynomial roots.

We remark that the integrand in (\ref{f:DL1}) can be approximated using the discrete difference in place of the derivative term:
\begin{equation}
z^p \frac{f'(z)}{f(z)}dz\rightarrow z_k^p\frac{f'(z_k)\Delta z}{f(z_k)}\approx z_k^p \frac{f(z_k)-f(z_{k-1})}{f(z_k)}=z_k^p \left(1-\frac{f(z_{k-1})}{f(z_k)}\right).
\label{f:logderapp}
\end{equation}
We will refer to the contour integration method used for the search of solitonic eigenvalues with integrand approximated as in (\ref{f:logderapp}) (i.e. without explicit derivatives) as aDL.

More recently, Kravanja et al in \cite{ksb99} presented an improved version of the DL approach. Their algorithm relied on a recursive construction of the so-called formal orthogonal polynomials, which have the roots that coincide with the zeros of our function $f(z)$. We have analyzed the accuracy of this improved method for eigenvalues computation too, but we do not present these results here on separate plots because our study did not reveal any noticeable difference between results obtained from this newer method compared to those of the ordinary DL method (at least, for the set of our test profiles). However, we note that for some real-world applications, where signals are not smooth and often significantly corrupted by noise \cite{tpl17}, the approach proposed in \cite{ksb99} might demonstrate a better performance.

Since all zeros inside the given contour can be located simultaneously together with their multiplicity, the contour integration method can guarantee more stability of the overall algorithm, especially when multiple eigenvalues are to be found. To ensure that we successfully localize all desired zeros, a large enough contour needs to be defined in the $\xi$ half-plane. In addition, the contour integration method runtime depends only insignificantly on the number of eigenvalues insofar as the values of logarithmic derivatives along the contour from (\ref{f:DL1}) can be computed just once and then saved for further processing. On the contrary, iterative algorithms, Subs.~\ref{subsec:iter}, can uncontrollably scan the complex plane and arrive at a zero value in the lower half-plane of $\xi$. In the case of several solitons, iterative algorithms evaluate all zeros independently, and that fact can increase the computational time. However, we note that for some well-defined cases, the contour integrals' computation can turn out to be more time-consuming in comparison with the iterative algorithms, because the former requires performing the computation of the value of $a(\xi)$ (possibly $a'(\xi)$ as well) for the whole set of discrete points along the contour $C$, while for a good initial guess the iterations may only require several repetitive evaluations of the function and its derivative.

\begin{figure}[h!]
	\centering
	\includegraphics[width=.6\textwidth]{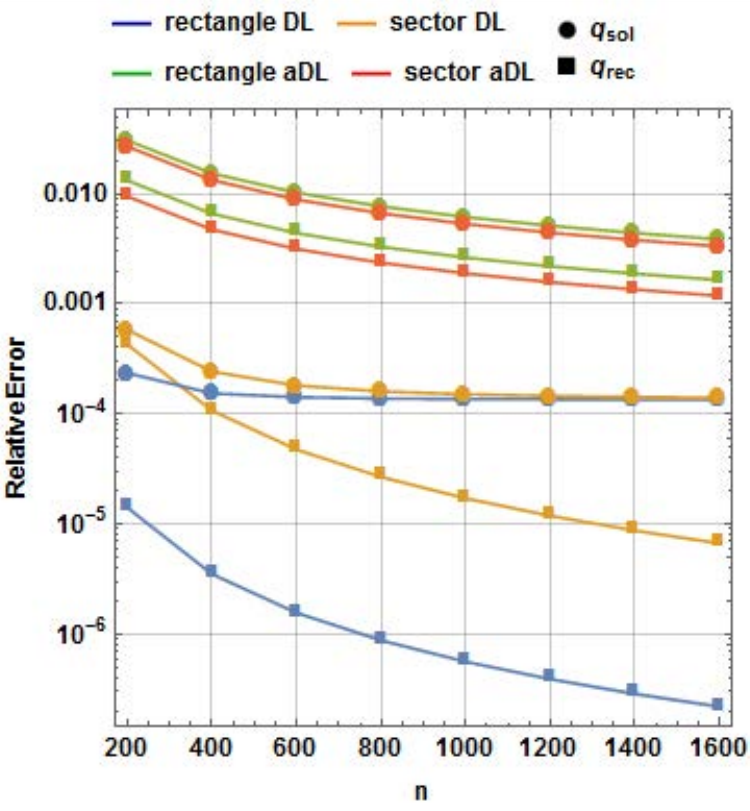}
	\caption{The dependence of the relative error on the number of contour discretization points
		for the rectangular potential ($L=1$, {\color{OrangeRed}{$A=\pi/2$  implies a single eigenvalue}}, \magenta{$2^{10}$ discretization points}) and for the soliton potential (\ref{f:shiftedpot})
		with phase factor ($L=20$, \magenta{$2^{10}$ discretization points}), for the rectangular integration contour: $\text{Re} (\xi)$ changes in the range $[-1 \ldots 1]$, $\text{Im} (\xi)$ changes in the range $[0.1 \ldots 2]$; and ring sector shape of the contour: $\xi =\rho e^{i \theta}$, $\rho$ changes in the range $[0.1 \ldots 2]$, $\theta$ changes in the range $[\pi/12 \ldots 11\pi/12]$. \magenta{For all curves BO method was used for ZSS solution.}}
	\label{fig:countour}
\end{figure}

We analyzed the dependence of contour integration algorithms' accuracy on the number of points along the contour and on the contour shape, Fig.~\ref{fig:countour}. We study the behavior of the aDL method with the approximation scheme represented in Eq. (\ref{f:logderapp}) and compare it with the ordinary DL method's behavior. As it was expected, the latter works more accurately due to a more accurate calculation of the derivative, see the blue and orange lines and compare them with red and green ones in Fig.~\ref{fig:countour}. We also observed that the particular contour shape chosen in Eq.~(\ref{f:DL1}) also influences the resulting accuracy of the eigenvalues found. In particular, we checked the behavior of the methods using the rectangular contour in the upper half-plane of $\xi$, fixing the contour borders along $\text{Re} (\xi)$-axis and $\text{Im} (\xi)$-axis. Another option that we tested was to define the ring sector in the $\xi$-plane, fixing the borders for absolute value $\rho$ and for argument $\theta$ of $\xi$ written in polar representation as $\xi = \rho e^{i \theta}$.

The integrals in all of our methods have been evaluated using the trapezoidal rule. Since the runtime of this algorithm depends linearly on the number of points along the contour, we found sufficiently enough to compare the runtimes of all approaches for the largest number of points (see Table~\ref{t:contours}). We found that results from the ordinary DL method and its approximated aDL version differed: the runtime of the ordinary DL method is typically 1.5-2 times larger than that of the aDL one, but the ordinary DL method allows us to reach a smaller relative error. From Fig.~\ref{fig:countour}, we can also readily see that the rectangular contour gives a better accuracy, at least for the test profiles that we studied. In our tests the computation time of the contour approach was typically higher than that for iterative algorithms from Subs.~\ref{subsec:iter} (cf. the runtimes for different methods given in Figs.~\ref{fig:poolrec},~\ref{fig:poolsol} and \ref{fig:countour}). This result is well explainable, since we were taking the guess points that were close enough to the zero point, such that the iterative algorithm was able to reach a zero in just several iterative steps, ensuring a lower computational time.

{\color{OrangeRed}{\begin{table}
	\caption{Runtimes (in seconds) of discrete eigenvalues' evaluation for different test profiles, contour integral approaches and shapes of contours: the rectangular integration contour: $\text{Re} (\xi)$ changes in the range $[-1 \ldots 1]$, $\text{Im} (\xi)$ changes in the range $[0.1 \ldots 2]$; and  ring sector shape of the contour: $\xi =\rho e^{i \theta}$, $\rho$ changes in the range $[0.1 \ldots 2]$, $\theta$ changes in the range $[\pi/12 \ldots 11\pi/12]$ for 1600 points along the integration contour.  \magenta{For all curves BO method was used for ZSS solution.}}\label{t:contours}
	\centering
	\def\arraystretch{1.2}
	\begin{tabular}{|c|c|c|c|}
		\hline
		rectangle DL&rectangle aDL&sector DL&sector aDL\\
		\hline
		\multicolumn{4}{|c|}{$q_{\text{sol}}$, \magenta{$2^{10}$ discretization points}}\\
		\hline
		88.27&51.49&88.27&51.37\\
		\hline
		\multicolumn{4}{|c|}{$q_{\text{rec}}$, $L=1$, $A=\pi/2$, \magenta{$2^{10}$ discretization points}}\\
		\hline
		109.65&69.75&109.65&69.64\\
		\hline
	\end{tabular}
\end{table}
}}.

\subsection{Derivative computation}\label{subsec:deriv}
Both the contour integration and the iterative algorithms (in particular, the NR method) require us to find the value of $a'(\xi)$ together with the value of $a(\xi)$ at the same point $\xi$. Since $a(\xi)$ is homomorphic in the upper half-plane of $\xi$ (see, for example, \cite{as}), the divided difference can be used to approximate the derivative. However, it is also possible to find the value of the derivative more accurately in the same programming loop together with the function computation itself. A similar approach is described in
\cite{hk16,yk14-2,bo92,bct98}. The idea here is to evolve the  derivative of $\Phi(t, \xi)$ along $t$ together with $\Phi(t, \xi)$ itself. In the case of ordinary ZSS (\ref{f:ZSode}), the evolution starts with the \enquote{initial condition} for the derivative defined at $t \to -\infty$: $\Phi'(t \to-\infty, \xi)\to(-it e^{-i\xi t},0)^T$. After the truncation of the $t$-interval, the initial condition takes the form $\Phi'(-L, \xi)=(iLe^{i\xi L},0)^T$. For the ZSS written for the envelope function (\ref{f:ZSodesimplified}), we have $X'(t\to-\infty, \xi)=(0, 0)^T$, or, after the truncation, $X'(-L, \xi)=(0 ,0)^T$. Evolution of the wave function derivative over a single step is performed by applying the relation:
\begin{equation}\label{der}
\Phi'_{m+1}=T_m'\Phi_m+T_m\Phi'_m,  \quad \text{or} \quad X'_{m+1}=T_m' X_m+T_m X'_m.
\end{equation}
Thus, for all transfer matrix methods from Subs.~\ref{subsec:meth}, the matrix $T'_m$ can be easily determined. At the same time, the RK method does not have that advantage: here we have to numerically solve the ZSS for the derivatives separately, equipped with different initial conditions.  However, we recall that the RK method is the most time-consuming according to our study, see Subs.~\ref{subsec:cont}, and so we do not deal with it further.

Our algorithm requires the use of the derivatives of the transfer matrices of each NFT algorithm described in Subs.~\ref{subsec:meth}. In the case of the BO matrix, the calculation of the derivative yields
\begin{equation}
T_m^{'\text{(BO)}}=\left(	\begin{matrix}A_+& B\\C&A_-\end{matrix}	\right),
\label{f:BOTprime}
\end{equation}
where
\begin{equation}
\begin{aligned}
&A_{\pm}=\pm\frac{i\xi^2\Delta t}{\kappa^2}\cosh \kappa \Delta t \mp  \sinh \kappa \Delta t \left(\pm\frac{\xi \Delta t}{\kappa}+\frac{i}{\kappa}+\frac{i\xi^2}{\kappa^3}\right),\\
&B=\frac{q_n\xi}{\kappa^3} \sinh \kappa \Delta t-\frac{q_n\xi\Delta t}{\kappa^2}\cosh \kappa \Delta t,\\
&C=-\frac{q_n^*\xi}{\kappa^3} \sinh \kappa \Delta t+\frac{q_n^*\xi\Delta t}{\kappa^2}\cosh \kappa \Delta t,
\end{aligned}
\label{f:BOTprime1}
\end{equation} where we keep the notations from Eq. (\ref{f:BOT}): $\kappa=\sqrt{-|q_m|^2-\xi^2}$. For the AL method, the derivative of the transfer matrix acquires the following form:
\begin{equation}
T_m^{'\text{(AL)}}=\frac{1}{\sqrt{1+\Delta t^2 q_m^2}}\left(
\begin{matrix}	
-i\Delta t e^{-i\xi \Delta t} &0\\
0 & i\Delta t e^{i\xi \Delta t}\end{matrix}	\right).
\end{equation}
For the modified BO method we have
\begin{equation}
T_m^{'\text{(BOmod)}}=2it\sin|q_m\Delta t|\left(\begin{matrix}0 &e^{i(\theta_{q_m}+2\xi t)}\\ e^{-i(\theta_{q_m}+2\xi t)}&0\end{matrix}\right),
\label{f:BOmodTprime}
\end{equation}
where $\theta_{q_m}$ is an argument of $q_m$. For the modified AL matrix we have:
\begin{equation}
T_m^{'\text{(ALmod)}}=\frac{2it\Delta t}{\sqrt{1+\Delta t^2 |q_m|^2}}\left(
\begin{matrix}	
0 & q_m e^{2i\xi t}\\
q_m^* e^{-2i \xi t} & 0
\end{matrix}\right).
\label{f:ALmodTprime}
\end{equation}
The CN-method transfer matrix derivative has the form \cite{yk14-2}:
\begin{equation}
T_m^{'\text{(CN)}}=\frac{1}{2}P_{m+1}(I-\frac{\Delta t}{2}P_{m+1})^{-2}(I+\frac{\Delta t}{2}P_m)+\frac{1}{2}(I-\frac{\Delta t}{2}P_{m+1})^{-1}P_m,
\label{f:CNTprime}
\end{equation}
where we again keep the notations from (\ref{f:CNT1}). In order to find $a'(\xi)$, an additional step must be taken at the end of the algorithm. For the ordinary ZSS (\ref{f:ZSode}), we have the following expression for the derivative of $a(\xi)$ that involves the elements of the Jost solution and its derivatives: $a'(\xi)=(\phi'_1(L, \xi)+i L\phi_1(L, \xi))e^{i\xi L}$; in the case of the envelope ZSS (\ref{f:ZSodesimplified}), we arrive at: $a'(\xi)=\chi'_1(L, \xi)$.

\begin{figure}[h!]
	\begin{tabular}{cc}
		\includegraphics[width=0.48\textwidth]{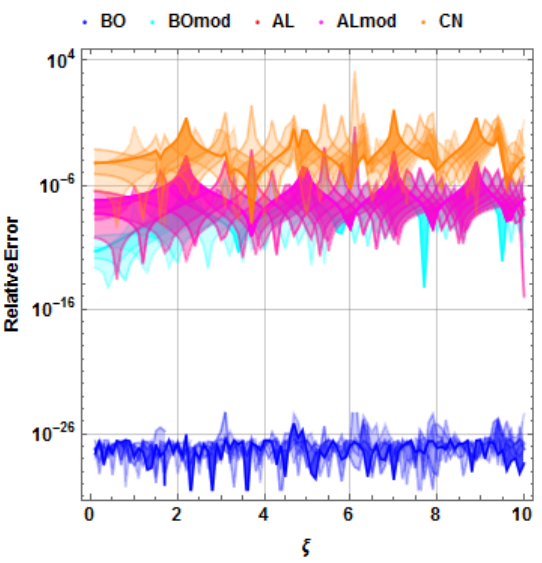}&
		\includegraphics[width=0.48\textwidth]{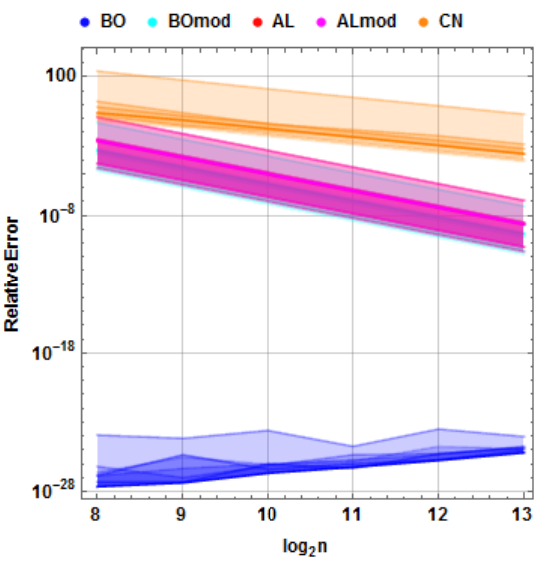}\\
		{\scriptsize{a)}}&	{\scriptsize{b)}}\\		
		\includegraphics[width=0.48\textwidth]{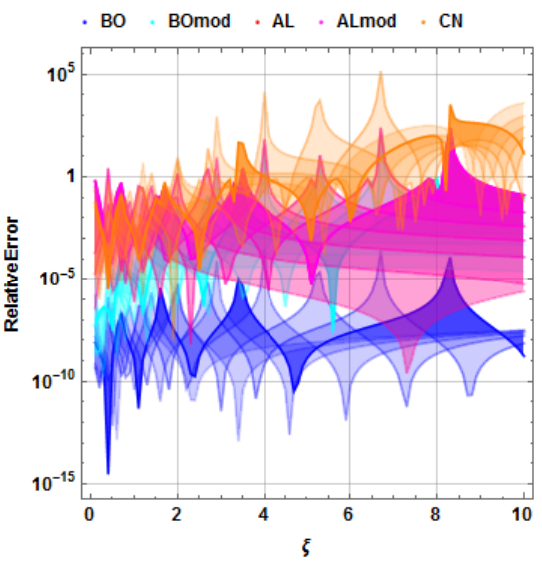}&
		\includegraphics[width=0.48\textwidth]{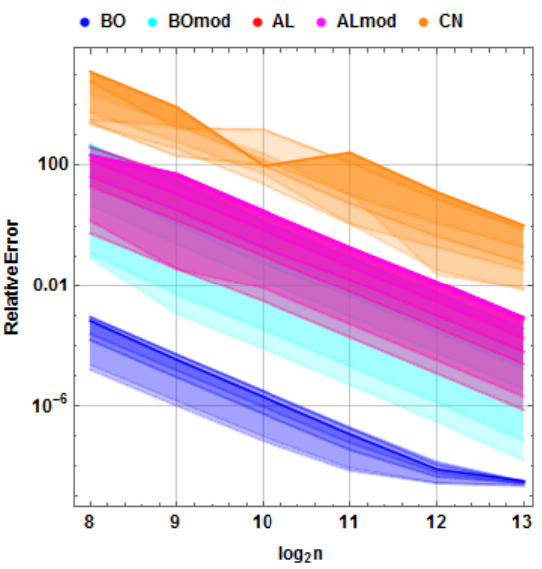}\\
		{\scriptsize{c)}}&	{\scriptsize{d)}}\\
		
	\end{tabular}
	\caption{Relative error  for $a'(\xi)$ computation versus $\xi$ {\color{OrangeRed}{(with step size $\Delta\xi=0.1$, $n=2^{10}$)}} and $n$ using different NFT algorithms from Subs.~ \ref{subsec:deriv}, panels a)--b): for the rectangular potential with $L=1$, amplitude changes in the range $[2, \, \ldots , \, 5]$ with the step $\Delta A=0.5$; panels c)--d): for the over-soliton potential with $L=20$, amplitude changes in the range $[2.25, \, \ldots, \, 5.25]$ with the step of the increment $\Delta A=0.5$.}
	\label{fig:aprime}
\end{figure}

To check the accuracy of the derivatives evaluated via different algorithms, we resort to similar mechanisms as described above. We primarily use the scattering coefficient $a(\xi)$ for the rectangular (\ref{f:reca}) and over-soliton (\ref{f:sola}) potentials. In Figs.~\ref{fig:aprime}, we present the dependence of the relative error for the derivative $a'(\xi)$ on $\xi$ and on the number of subintervals $n$ for the different methods described above. Our computations confirms that the BO method is again the most accurate one, while the accuracy of the CN is the lowest among all methods studied.
\subsection{Computation of the total number of eigenvalues}\label{subsec:zeros}
A serious difficulty in dealing with iterative algorithms is that they require an additional adaptation for the multisolitonic case. By definition, the iterative algorithms seek for only one zero of $f(x)$, which is usually the closest one with respect to the initial guess point. If the input signal comprises of more than one eigenvalue (e.g. when we have sufficiently large amplitudes for the rectangle or over-soliton potentials), it leads to the following requirements for iterative algorithms:
\begin{itemize}
	\item The total number of eigenvalues must be known before the search routine is executed. The number of eigenvalues can either be known from the system properties or, alternatively, must be computed.
	\item When a solution (zero value of the function) is found, it must be eliminated from the next search runs, to exclude the possibility of repetition.
\end{itemize}

\begin{figure}[h!]
	\centering
	\includegraphics[width=0.75\textwidth]{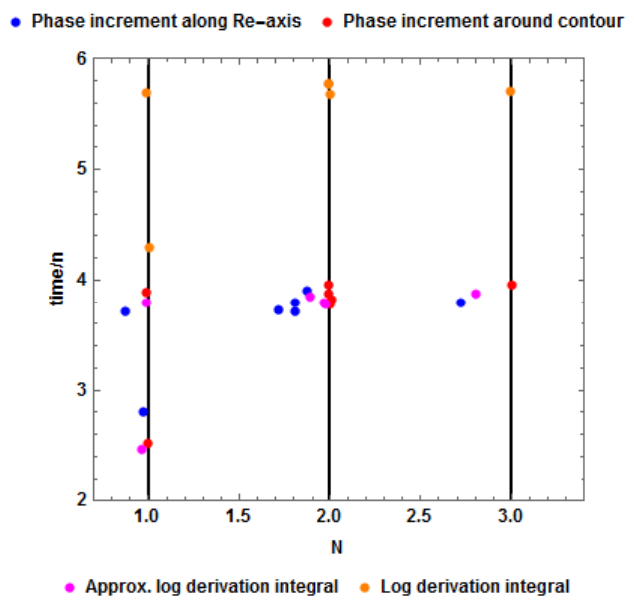}
	\caption{ Effective ``number of zeros'', obtained via the computation of phase increment, versus runtime diagram for different number of zeros evaluation techniques (the BO NFT computation method was used for the illustration) {\color{OrangeRed}{at the following profiles: the soliton with phase factor $q_{\text{sol}}$, over-soliton $q_{\text{over}}$ and rectangular potentials, both with $A=2$, $2.5$ and $3$}}.\label{fig:nzeros}}
\end{figure}

In order to find the total number of zeros, we can use the following relation (logarithmic derivative):
\begin{equation}
N=\frac{1}{2\pi i}\int_C \frac{f'(x)}{f(x)}dx=\frac{1}{2\pi}\Delta_C \text{arg}f(x).
\label{f:N}
\end{equation}
In Eq. (\ref{f:N}) $\Delta_C \text{arg}f(x)$ represents the incremental change in the value of the argument of a complex-valued function $f(x)$, when $x$ traverses the closed contour $C$ in the complex plane. We can manually choose a sufficiently large $C$ in order to keep all possible zeros inside the contour. Alternatively, we can calculate phase increment along the real axis, assuming that the contour contains a part of the real axis and that the other part of the contour located  in the upper half-plane makes negligible contribution to the total increment value. We tested  the applicability of all described ways for the computation of the number $N$ of zeros and found out (see Fig.~\ref{fig:nzeros}) that the precise logarithmic derivative integration was the slowest option for determining $N$.
 Then, the computation of the phase increment and the approximate integration using Eq.~(\ref{f:logderapp}) as integrand in Eq.~(\ref{f:N}), both work faster and produce a smaller deviation from an integer number; the latter are displayed in Fig.~\ref{fig:nzeros} as thick vertical lines. To plot this diagram, we check the number of zeros for all presented model signals from Subs.~\ref{subsec:signal} with randomly chosen amplitudes. We notice that for the computation of the total number of eigenvalues, $N$, it is not necessary to compute the expression (\ref{f:N}) very accurately; the numerical error just needs to be small enough to distinguish between two successive integer values of $N$.

In order to eliminate an estimated zero, $x_i$, from the next search run in the iterative schemes, the function $f(x)$ has to be redefined as
\begin{equation}\label{f:elim}
f(x)\rightarrow \frac{f(x)}{x-x_i}.
\end{equation}
We note that after the redefinition of $f(x)$, there can arise some convergence problems, if the search path goes in the vicinity of already located zeros. In addition, the approach cannot help us to eliminate multiple zeros, they will be found according to their multiplicity. The latter could be useful only if we need to know the eigenvalues multiplicity, otherwise the algorithms will waste time to find an already found zero. We remark that for the real optical applications, it is unlikely to have multiple zeros of the scattering  function $a(\xi)$ due to the presence of noise that should typically split the multiple eigenvalue. Because of that issue we do not pay much attention to the multiplicity in the computation of  eigenvalues in our present study.

At this point, it is pertinent to make a remark on the region, in which the eigenvalues are sought for. The values of real and imaginary parts of an eigenvalue can be limited by physical and mathematical reasons. By definition, the  eigenvalues lie in upper half-plane, so $\text{Im} (\xi)>0$ gives one geometrical bound for the search region. At the same time, the imaginary part corresponds to the amplitude of the soliton, which is associated with its energy.
Hence, a knowledge of the energy of the signal allows us to impose an upper limit on the boundary of $\text{Im} (\xi)$. We note that for the contour integration algorithms we seek the eigenvalues inside the manually defined region of interest (ROI) in the complex $\xi$-plane. For iterative algorithms, the iteration paths are not usually controlled (unless some constraints inside the search routine are additionally imposed) and thus the iteration algorithm scan incidentally go out of the ROI; this can be an additional drawback for that group of methods.

\subsection{Multisolitonic test}\label{subsec:mult}
Multisolitonic potentials are patently interesting from the perspective of practical applications, with subsets of multisolitonic eigenvalues being specifically proposed for optical communication purposes \cite{tm13,hn93,hky14,hk16}. In this subsection, we investigate the performance of both iterative and contour integration methods for the case when the NFT pulse decomposition involves several eigenvalues, and we need to retrieve their values. For our tests, we choose the over-soliton potential as in Eq. (\ref{f:solpot}) for the fixed amplitude $A=5$. This signal has five eigenvalues in its discrete NF spectrum:
\begin{equation}\label{f:eig}
\xi_k=(4.5-k)i, \quad k=0 \ldots 4.
\end{equation}
The full decomposition also contains the non-zero continuous spectral data, so that the situation considered in this section is quite general. To locate the eigenvalues, it is convenient to use the rectangular ROI for the contour integration methods because, as it was found in Subs.~\ref{subsec:contour}, such a contour ensures a greater accuracy in results. We choose the contour having the shape of a rectangle in the complex $\xi$-plane with the dimensions: $[-1,1]$ along the real axis and $[0.1,5]$ along the imaginary axis. Such a rectangle encompasses all eigenvalues defined by Eq.~(\ref{f:eig}) for our  over-soliton profile.

Now we note that our iterative algorithms need adaptation for the multisolitonic search task, as described in the previous subsection. The adapted multisoliton-search iterative algorithm scheme is given in Alg.~\ref{alg}.
\begin{algorithm}
\caption{Adaptation of an iterative algorithm to the multisolitonic case}\label{alg}
\begin{tabular}{l}
	\hline
	{\bfseries{Input:}} Define function $f(x)$, make initial guess $x_0$, \\expected number of zeros $N$, ROI.\\
	{\bfseries{Step 0:}} Initialize an empty output array $x_{out}$.\\
	{\bfseries{Step 1:}} Define a current function $f_{c}(x)=f(x)$ \\and a current guess $x_{c}=x_0$.\\
	{\bfseries{Step 2:}} Launch the iterative algorithm for $f_{c}(x)$, $x_{c}$\\ and limited number of iteration steps.\\
	{\bfseries{Step 3:}} Check if zero $x_i$ was successfully located in the previous step. \\
	If yes, go to step 4. \\
	If no, go to step 7.\\
	{\bfseries{Step 4:}} Check if the located zero is inside the ROI. \\
	If yes, concatenate the located zero to the output array $x_{out}=[x_{out}, x_i]$ \\and go to step 5.\\
	If no, go to step 7.\\
	{\bfseries{Step 5:}} Check if all expected zeros are located using {\texttt{size($x_{out}$)}}==$N$.\\ If yes, go to Output.\\If no, go to step 6.\\
	{\bfseries{Step 6:}} Redefine the  function to eliminate located zero\\ $f_{c}(x)=f_{c}(x)/(x-x_i)$.\\
	{\bfseries{Step 7:}} Choose an initial guess number $x_{c}$ randomly from inside the ROI \\and go to step 2.\\
	{\bfseries{Output:}} Estimate the array of zeros $x_{out}$.\\
	\hline
\end{tabular}
\end{algorithm}

\begin{figure}[h!]
		\centering
		\includegraphics[width=0.5\textwidth]{Figlegend}
			\includegraphics[width=0.65\textwidth]{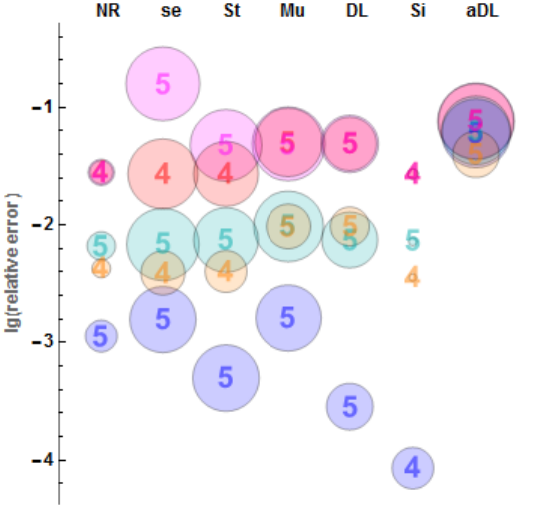}
	\caption{\color{OrangeRed}The bubble chart shows the relative error of the eigenvalue location computed via the different root-finding iterative (NR, se, St, Mu, Si) and contour integral (DL, aDL) methods for over-soliton profile with five embedded solitonic eigenvalues as in Eq.~(\ref{f:eig}). The digit in each bubble means the number of ultimately located zeros (the maximum is 5), the error was assessed as the mean relative error via Eq.~(\ref{f:relerror}), the formula was applied for located zeros only. The size of each particular bubble relatively shows the inverse runtime, a smaller bubble means a longer run and vice versa; the logic applies to concentric and overlapping circles. Runtime changes from 0.025 s (the largest bubble) to 4192 s (the smallest one). The color of each bubble identifies to which NFT method the root-finding algorithm was coupled.}
	\label{fig:multisoliton}
\end{figure}

In our specific case, we started with the initial guess $\xi_0=i$, followed by up to three attempts to locate the same zero over $10^5$ iteration steps to reach each particular zero point.
The resultant error versus runtime diagram is represented as a bubble chart in Fig.~\ref{fig:multisoliton} for the different combinations of the transfer-matrix NFT algorithms from Subs.~\ref{subsec:meth} and the root-finding iterative methods from Subs.~\ref{subsec:iter}~and~\ref{subsec:contour}. Here bubble sizes are inversely proportional to the runtime with the numbers therein indicating the number of zeros; respective colors identify the numerical algorithm chosen.

We see that, in accordance with the remarks made in the previous subsection, the iterative algorithms equipped with the elimination procedure as in Eq.~(\ref{f:elim}) may fail to estimate the entire cluster of zeros even when coupled with the most accurate BO NFT method. At the same time, the contour integration methods, especially aDL, give a higher error margin.

\subsection{A hybrid method}\label{subsec:hybr}
As we noted in the previous parts of this section, both types of eigenvalue-finding approaches have some inherent disadvantages. The iterative algorithms from Subs.~\ref{subsec:iter} can be unstable: we cannot be confident that the iterations will eventually lead to the correct eigenvalue points (see Figs.~\ref{fig:poolsol}~and~\ref{fig:poolrec}). They also require some additional adaptation to the multisolitonic case to incorporate the elimination of previously found zeros, and the restriction on the search region has to be generally imposed. The contour integration algorithms from Subs.~\ref{subsec:contour} do not allow us to reach a high accuracy and take a comparatively long time.

 In this subsection, we present the new hybrid method, which allows us to take advantage of the best from both approaches whilst simultaneously getting rid of their respective drawbacks. The main idea of our hybrid method is that \textit{we can use the result of contour integration as the initial guess that is then supplied to the consequential iterative algorithm} {\color{OrangeRed}(the same strategy was mentioned in \cite{dl67} as a way of root-finding refining)}. This combination allows us to reach almost any accuracy up to the limitations imposed by the NFT computation method itself, see Subs.~\ref{subsec:meth}.
{\color{OrangeRed}{{The hybrid method presented here guarantees locations for all eigenvalues as opposed to iterative methods that often found most eigenvalues but not necessarily all.}}} The description of the consecutive steps for the hybrid algorithm is as follows.
\begin{itemize}
\item First, find the location of the approximate zeros' using one of the contour integration method drawing a large enough contour. Two key remarks on this step here:
\begin{itemize}
\item the integration result does not require to be really accurate, so the computation time can be reduced;
\item it allows us to find a good approximation for all zeros that lie inside our ROI;
\end{itemize}
\item Second, we apply a particular iterative method to find more precise location for each eigenvalue, using the results of the previous step as guess points (and, eventually, also employing some other data that can be obtained by the contour integration, i.e. the multiplicity).
\end{itemize}

\begin{figure}[h!]
\centering
\includegraphics[width=0.4\textwidth]{Figlegend}
\begin{tabular}{cc}
	\includegraphics[width=0.48\textwidth]{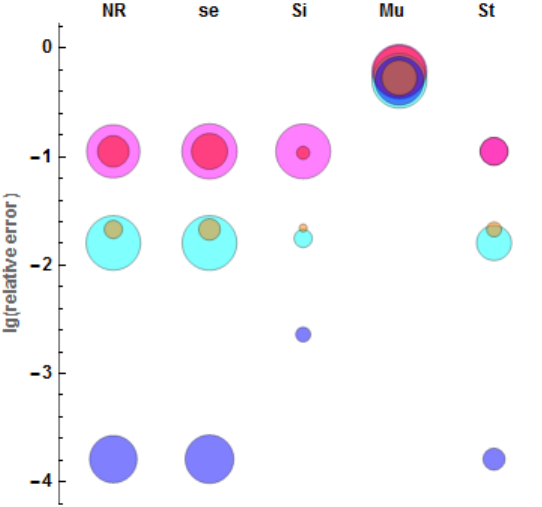}&	\includegraphics[width=0.48\textwidth]{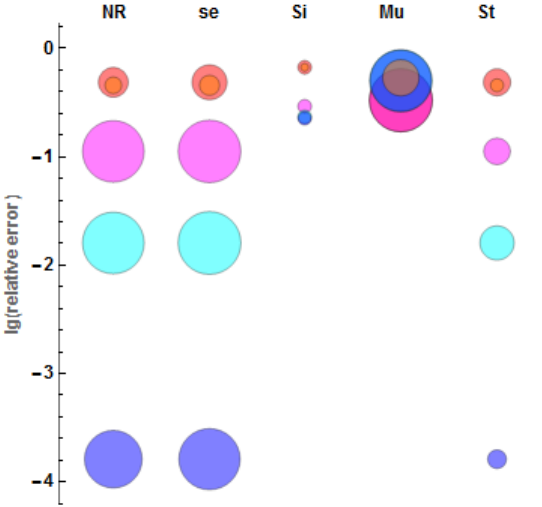}\\
	{\scriptsize{a) DL}}&	{\scriptsize{b) aDL}}
\end{tabular}
	\caption{ {\color{OrangeRed}{The bubble chart shows the relative error of solitonic eigenvalue location by using the hybrid method; the results are given for the over-soliton profile with five eigenvalues. Panel a) corresponds to the DL contour integration to find the guess values; panel b) does the same by using the aDL method. As before, the bubble size shows inverse runtimes, a smaller bubble indicates a longer run and vice versa, the logic applies to both concentric circles. Runtime changes from 3.27 s (the largest bubble) to 286 s (the smallest one). The color of each bubble identifies to which NFT method the root-finding algorithm was coupled.}}}
	\label{fig:hybrid}
\end{figure}

To test our new hybrid algorithm, we again employ the over-soliton potential from the previous subsection \ref{subsec:mult} with five solitonic modes, and plot the error of the eigenvalues computed in Fig.~\ref{fig:hybrid}. We test the different combinations of the NFT algorithms from Subs.~\ref{subsec:meth} (we do not use the RK method in this section as it is too time consuming) combined with contour integration methods from Subs.~\ref{subsec:contour} and iterative methods from Subs.~\ref{subsec:iter}.

\begin{table}
	\magenta{\caption {Normalised runtimes (in milliseconds per sample of the signal) of solitonic eigenvalues evaluation for multisoliton profile ($A=5.25$, $n=2^{10}$) for iterative and contour integral methods, implemented individually (first row and first column) and for hybrid algorithm, when these methods are paired. BO method was used for ZSS solution. Rectangular contour with 1600 discretization points was used for contour integration.\label{t:multitimes}}}
	\centering
	\def\arraystretch{1.2}
\magenta{\begin{tabular}{|c||c|c|c|c|c|}
		\hline
		~&\textbf{NR}&\textbf{secant}&\textbf{Steffensen}&\textbf{Muller}&\textbf{Sidi}\\
		~&19.18&3.56&3.53&3.56&10.22\\
		\hhline{======}
		\textbf{DL}&7.68& 7.3& 76.83&7.31&35.33\\
		7.8&~&~&~&~&~\\
		\hline
		\textbf{aDL}& 3.88&3.44&71.75&3.41&35.86\\
		3.28&~&~&~&~&~\\
		\hline
	\end{tabular}}
\end{table}
{\color{Black}{We see that, in general, runtime reduces comparably with contour integrals approach, but it is still larger than runtime of some iterative algorithms (see Table~\ref{t:multitimes} for the runtimes of BO method used as example).}} We see that the largest runtime is demonstrated by the Sidi and Steffensen methods, but the latter show a surprisingly high accuracy when coupled with the BO NFT algorithm. The NR, Muller and secant algorithms perform fairly similar in terms of runtime, but the Muller method displays a worse accuracy. If we combine the NR and secant methods with any of the integral methods, then the utilization of both the DL and aDL approaches gives similar results. It means, that both the DL and aDL methods supply a sufficiently good initial approximation to reach our desired zero.

\section{Computation of norming constants}\label{sec:norm}
The third and last component of the NF spectrum is the norming constants (\ref{f:res}) attributed to each eigenvalue. These parameters define the phase and the center position of each solitonic degree of freedom \cite{as,zs72,yk14-1}. Under the assumptions listed in Secs. \ref{sec:intro},  \ref{sec:prem}, the norming constants are expressed as the residues of the reflection coefficient $r(\xi)$ from (\ref{f:r}) calculated at its simple poles, Eq.~(\ref{f:res}), i.e. at the solitonic eigenvalues addressed in the previous section. So at this point we assume that some appropriate method from Sec.~\ref{sec:disc} has been executed and the plausible values for all eigenvalues are now known with high enough precision.

\subsection{Ways of the residues evaluation}\label{subsec:norm-mrth}
For the computation of the norming constants, we can equally use each of the two expressions for the residue. First, the norming constant can be found via the contour integration:
\begin{equation}
	c_j=\frac{1}{2\pi i}\int_{\gamma_j}\!\!r(\xi) \, d\xi,
		\label{f:resintegral}
\end{equation}
where $\gamma_j$ is a sufficiently small contour in the $\xi$-plane encircling the single pole $\xi_j$ (considered to have been already located). Alternatively \cite{tpl17,yk14-1}, the norming constant is given by fraction
\begin{equation}
	c_j=\frac{b(\xi_j)}{a'(\xi_j)},
	\label{f:resfraction}
\end{equation}
where $a'(\xi_j)$ is derivative of the scattering coefficient $a(\xi)$ from (\ref{f:ab}) with respect to $\xi$ evaluated in its simple zero $\xi_j$.  The computation of $a(\xi)$, $b(\xi)$ or $r(\xi)$ in any point of $\overline{{\mathbb C}^+}$ can be performed using the methods from Subs.~\ref{subsec:meth}. The computation of $a'(\xi)$ can be carried out in the same program cycle in parallel with the other NFT quantities that do not involve derivatives as it is described in Subs.~\ref{subsec:deriv}. However, as we shall see, the computation of the norming constants brings about some additional problems \cite{hk16,a16}.

\subsection{Test of the norming constants straightforward computation}\label{subsec:norm-test}
Explicit analytical expressions for the norming constants corresponding to the eigenvalues of three model potentials that we use in this paper are presented in Subs.~\ref{subsec:signal}. For the forthcoming tests we use the over-soliton profile (\ref{f:solpot}), for which the expressions for both eigenvalues and norming constants are explicit and relatively simple. As for the rectangular potential (\ref{f:recpot}), the eigenvalues are defined through the solutions of a transcendental equation, which, in turn, can be found only with a finite accuracy, and this can restrict our comparative analysis. Therefore we do not use the rectangular potential for testing the norming constant methods accuracy.

\begin{figure}[h!]
	\centering
	\includegraphics[width=0.4\textwidth]{Figlegend}
	\begin{tabular}{cc}
		\includegraphics[width=0.48\textwidth]{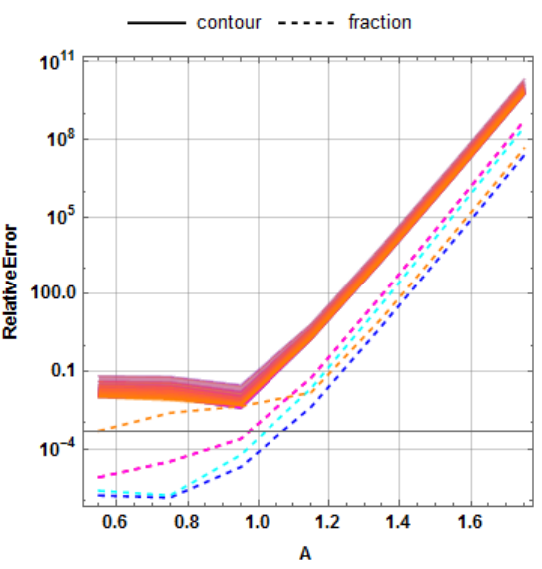}&
		\includegraphics[width=0.48\textwidth]{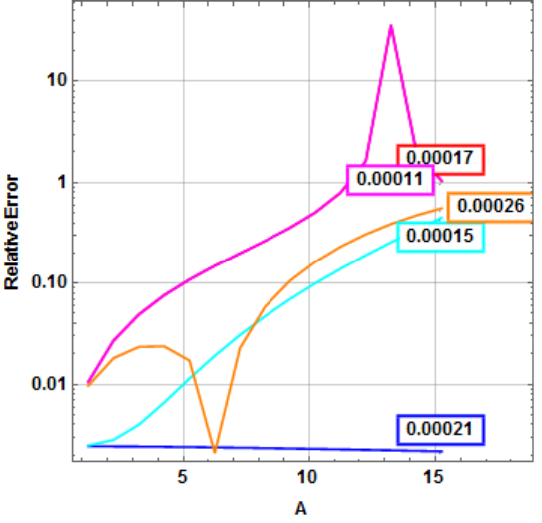}\\
		{\scriptsize{a)}}&	{\scriptsize{b)}}
		\end{tabular}
	\caption{Residue relative error on over-soliton signal's amplitude for different NFT algorithms and residues computation methods ($n=2^{12}$, $L=20$) for pane a) fraction and contour integral formulae (number of points along integration contour rise from transparent to solid lines in the range $n_{p}=[20\dots100]$ with the increment step $\Delta n_{p}=20$); panel b) fraction formula with application of improved scheme of $b(\xi)$ evaluation.}
	\label{fig:residues}
\end{figure}

In our test, we compare the performance of the different NFT transfer-matrix algorithms from Subs.~\ref{subsec:meth}~and~\ref{subsec:deriv}; these are combined with the two residue expressions given in the previous subsection. From Fig.~\ref{fig:residues}a) we can see that the fraction formula (\ref{f:resfraction})  employed for the residue computation, improves the result accuracy as compared to the integral definition (\ref{f:resintegral}) (both combined with different NFT methods). The main feature that we can extract from the analysis of Fig.~\ref{fig:residues} is that beyond a critical value of the amplitude ($A\lesssim1$ for this particular case), the computational error starts to rise rapidly. Manipulating the number of discretization points for the computation of the integral from~(\ref{f:resintegral}) does not change this error signature, but it evidently influences the runtime significantly. The integral formula (\ref{f:resintegral}) gives indistinguishably similar results independent of the NFT algorithm chosen. The fraction formula~(\ref{f:resfraction}) gives the best accuracy for the BO method and the worst accuracy attributed to the CN one. In general, our current results reveal that both formulae of the residue combined with different available NFT methods are relatively inaccurate, and this can be an important degrading factor in practical applications.

Authors of \cite{yk14-2} also showed that the AL, CN, Euler and central differences methods applied for the ZSS solution gave the similar unimpressively-small accuracy of norming constant computation. In \cite{bct98}, it was also found that here the RK method fares worse than the BO approach. In \cite{a16}, the author used the BO method and compared its performance with AL, CN for the norming constants computation.

\subsection{Discussion and the improvement of norming constants computation accuracy}\label{subsec:discus}
The problem of the numerical computation of residues lies in the properties of the spectral
functions $a(\xi)$, $b(\xi)$ and $r(\xi)$. When performing the analytical continuation of these functions into the complex $\xi$-plane, we should carefully check when this operation is indeed legitimate. As stated above,
the analytic continuation of  $b(\xi)$ is limited, in general, by the rate of decay of the ZSS potential.
On the other hand, when dealing with the function $a(\xi)$, its analytic continuation can be performed over the whole $\mathds{C}^+$ \cite{as},
where it is bounded and, moreover, $a(\xi)\to 1$ as $\xi\to\infty$, and so this property adds on our motivation to use the function $a(\xi)$ for the solitonic eigenvalues.

In the case when  $b(\xi)$ admits the analytic continuation into the whole complex plane
or a part of it covering the location of the eigenvalues $\{\xi_j\}$, 
the norming constants $b_j$ associated to $\xi_j$ can be calculated as (see
(\ref{f:eigen})) $b_j=b(\xi_j)$.
However, $b(\xi)$ is likely to grow exponentially when the imaginary part of $\xi$ increases. Numerical computations are also limited by the speed of the processor and by the available memory size (e.g. double precision numbers are bounded by approximately {\texttt{1.8\text{e}308}}). When we reach this limiting value, the accuracy of the computation would be naturally affected.  In order to understand the scale of problem arising in the accurate norming constant computation, we compared the limitations caused by the analytic continuation problems and by the computational reasons listed above.

\begin{figure}[h!]
	\centering
	\includegraphics[width=0.8\textwidth]{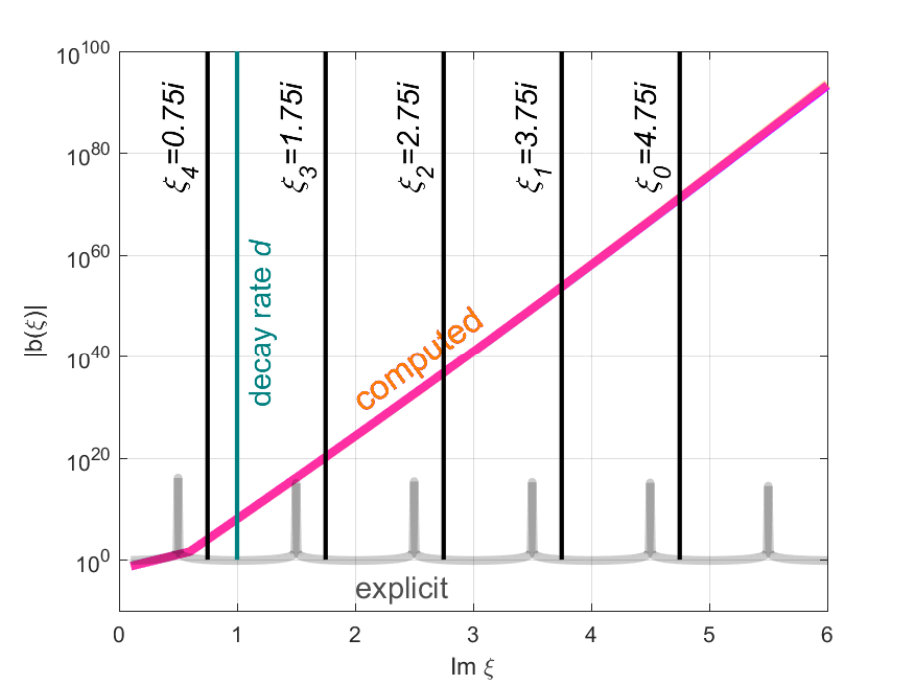}
	\caption{The computed values of $b(\xi)$ (red line) along imaginary axis of $\xi$, the value of the explicit expression (gray line), eigenvalues (black lines), and the potential decay rate (green line), calculated for the over-soliton profile (\ref{f:solpot}) with $A=5.25$.}
	\label{fig:borders}
\end{figure}

For the over-soliton potential with amplitude $A=5.25$, we calculate the values of $b(\xi)$ for a purely imaginary $\xi$. The curves for the numerically computed $b(\xi)$ (see  Fig.~\ref{fig:borders}) grow extremely fast, and so it is not possible to distinguish individual lines for the different NFT algorithms. It is also seen that the rise of $b(\xi)$ starts well before the decay rate $d$ coming to play (it is responsible for the band in the $\mathbb{C}$, where $b(\xi)$ is certainly analytical and shown by the green line). Therefore, our NFT algorithms are unable to evaluate the correct value of $b(\xi)$ (the gray line on the plot) even for sufficiently small $\text{Im}\xi$.  On the other hand, the evaluation of the norming constants 
in terms of $b(\xi)$
requires the values of $b(\xi)$ at the eigenvalues (marked by the black lines), where, as we see, the spectral function $b(\xi)$ is computed incorrectly.

However, we can adapt the NFT algorithms to provide a more accurate computation of $\{b_j\}$ 
associated with the eigenvalues $\{\xi_j\}$. Indeed, the detailed control of each computational step shows that if one uses Eq.~(\ref{f:absimpl}) for the calculation of $b$, the aforementioned divergence happens at the truncation edges of the potential (exponential factor grows fast for large $\text{Im}(\xi)$ and $L$).
Thus it is possible to simplify the evaluation of $b_j$ (i) by estimating it at the centre of the evaluation interval, avoiding the edges,
and (ii) by getting rid of the exponentials in the initial conditions as in 
Eq.~(\ref{f:asy}).  This idea was effectively harnessed in \cite{hk16} for the solution of ZZS employing a more straightforward Euler method and in \cite{a16} for BOmod transfer-matrix approach. It can be realised by considering the evolution of the so-called left wave (defined by its asymptotic as $t\to-\infty$) from the left edge towards the centre of the interval, and the evolution of the right wave (defined by its asymptotic as $t\to+\infty$) from the right edge to the centre towards the left wave,
and then by using Eq.~(\ref{f:eigen}) at $t=0$ to calculate $b_j$.

Using the wave function envelopes for both $\Phi(t, \xi)$ and $\Psi(t, \xi)$, defined as
\begin{equation}
G(t, \xi)=\Phi(t, \xi) e^{i\xi t}, \qquad
H(t, \xi)=\Psi(t, \xi) e^{-i\xi t},
\label{f:psieta}
\end{equation}
we can evaluate the coupled systems for both vectors $G(t)$ and
$H(t)$:
\begin{equation}
\frac{d}{dt}G(t, \xi)=\left(\begin{matrix}0&q(t)\\-q(t)^*&2 i\xi \end{matrix}\right)G(t, \xi),
\end{equation}
\begin{equation}
\frac{d}{dt}H(t)=\left(\begin{matrix}-2 i \xi&q(t)\\-q(t)^*&0 \end{matrix}\right)H(t, \xi).
\end{equation}
After the potential truncation, the initial value of the vector $H(L, \xi)=(0, 1)^T$ evolves from $t=L$ towards $t=0$, whilst the vector $G(-L, \xi)=(1, 0)^T$ evolves from $t=-L$ to $t=0$. At the point $t=0$, the desired quantity $b_j$ can be obtained from the following relation (cf. (\ref{f:eigen})):
\begin{equation}
G(t=0, \xi_j)=H(t=0, \xi_j)b_j.
\end{equation}
The evolution of $G(t, \xi)$ and $H(t, \xi)$ can be performed similarly to Subs.~\ref{subsec:meth}, i.e. using the transfer matrices. The implementation of the BO approach \big({\color{OrangeRed}{we used matrix exponential of ODE matrix,}} see Eqs.~(\ref{f:BOT1}), (\ref{f:BOT})\big) leads to the following transfer matrix for the left wave
\begin{equation}
T_m^{\text{(BOleft)}}=e^{i \xi \Delta t}\left(	\begin{matrix}	\cos \kappa \Delta t - i\xi / \kappa \sin\kappa \Delta t & q_n/\kappa \sin \kappa \Delta t\\
-q_n^*/\kappa \sin \kappa \Delta t & 	\cos \kappa \Delta t + i\xi / \kappa \sin \kappa \Delta t\end{matrix}	\right)
\end{equation}
and
\begin{equation}
T_m^{\text{(BOright)}}=e^{-i \xi \Delta t}\left(	\begin{matrix}	\cos \kappa \Delta t - i\xi / \kappa \sin\kappa \Delta t & q_n/\kappa \sin \kappa \Delta t\\
-q_n^*/\kappa \sin \kappa \Delta t & 	\cos \kappa \Delta t + i\xi / \kappa \sin \kappa \Delta t\end{matrix}	\right)
\end{equation}
for the right wave (here  $\kappa=\sqrt{|q_m|^2+\xi^2}$).
The AL-type approach \big({\color{OrangeRed}{(using Euler method and exponent first order decomposition,}} see~(\ref{f:ALT})\big) leads to the following transfer matrices for the left and right envelopes, correspondingly:
\begin{equation}
T_m^{\text{(ALleft)}}=\frac{1}{\sqrt{1+\Delta t^2 q_m^2}}\left(	\begin{matrix}	1 & q_m \Delta t\\
-q_m^*\Delta t & 	e^{2 i \xi \Delta t}\end{matrix}	\right),
\end{equation}
\begin{equation}
T_m^{\text{(ALright)}}=\frac{1}{\sqrt{1+\Delta t^2 q_m^2}}\left(	\begin{matrix}	e^{-2 i \xi \Delta t} & q_m \Delta t\\
-q_m^*\Delta t & 1	\end{matrix}	\right).
\end{equation}
We emphasize that the presented way of  definition and computation of 
$b$
is appropriate only for the discrete ZSS eigenvalues.

We tested our new methods of calculation of $\{b_j\}$ to find the residues through the fraction formula presented in Eq.~(\ref{f:resfraction}) (with $b(\xi_j)$ replaced by
$b_j$); our results are summarised in Fig.~\ref{fig:residues}b). We can readily observe that by employing the method described above we get more accurate norming constants for the wider range of the amplitude of the potential. In contrast to our previous results, the CN method shows better accuracy than both AL and ALmod methods, but the CN method takes more time than all other methods. The BO method gives the best accuracy and the weakest dependence of the outcome on the amplitude.

\section{Conclusion}\label{sec:concl}
In this study, we have formulated a comprehensive list of options for NFT computation, and then have made categorical comparisons of their relative advantages and disadvantages, principally focusing on runtime and accuracy optimization. The focus here has been to obtain best possible numerical approximation algorithm of the given spectral data.

First, the calibration of different transfer-matrix algorithms' performance was done using the continuous part of the spectral data as a measure for the method's accuracy analysis. In particular, it was shown that the BO method is usually superior to other alternatives in terms of the runtime and the accuracy of the result obtained, a conclusion that complies with some earlier existing studies on the NFT methods performance.

However, a major incentive of this work was the efficient computation of the  eigenvalues and norming constants associated with each solitonic degree of freedom (i.e. in the full \textit{discrete spectral data} associated with a given profile). For the computation of the eigenvalues, we first applied different iterative algorithms (involving the derivative computation or avoiding it at each step) combined with different transfer-matrix methods for the ZSS solution, using three different profiles with the known NF spectrum. At the beginning, we analyze the regions of convergence of each NFT method combined with a particular iterative scheme that relied on a reasonably accurate guess point that is then used for a search run. Surprisingly enough, we find that the Muller method typically provides the largest convergence region; the relatively large convergence basin is observed for the derivative-free Sidi method as well, while the Newton-Raphson, secant and Steffensen's iterative approaches show progressively poorer convergence.

This is followed by a description of a new class of methods for the search of eigenvalues
(the zeros of the spectral function $a(\xi)$) based on the evaluation of contour integrals in the complex plane of the spectral parameter $\xi$. Although we found that these methods are often slower and less accurate compared to \enquote{more traditional} iterative search algorithms, they generally allow us the freedom of a lax choice for the initial guess point during the eigenvalue search. This task is often difficult to fulfil in realistic applications, where the properties of the priories are not known a priori. In this section, we present a new hybrid algorithm that allows us to combine merits of both the iterative and contour integral method classes. The new method relies on the approximate initial evaluation of the NFT data using the contour integrals, and then using the obtained data as a good guess for an iterative algorithm following the initial search. This scheme leads to faster evaluation of the eigenvalues with high enough accuracy while largely being able to avoid convergence and initial guess issues. This new method for eigenvalues search is one of the main results of our present paper.

Finally, we address the problem of the computation of the norming constants using different algorithms. As same as the earlier studies suggest \cite{hk16,yk14-2}, the straightforward application of the single-directional NFT algorithms usually results in lower accuracy of the estimated norming constants. This adverse finding can be attributed to the properties of the analytic continuation of the function $b(\xi)$ into the upper half-plane of $\xi$, which is used in the norming constant evaluation. To rectify this problem, we utilise a numerical scheme first suggested in \cite{hk16} that employs simultaneous right and left scattering. We then use this bi-directional approach combined with different transfer-matrix methods. Our results prove that this approach results in an improved accuracy of the norming constants computation.

We believe that this new regime of modified NFT algorithms has led to a framework that will concomitantly optimize speed against accuracy, thereby to drive an adjustable numerical routine for the evaluation of NF spectral data tailored specifically to the problem in hand. This subjectivity aspect can be highly valuable in view of the active current progress of optical communication methods based on the NF spectrum modulation \cite{tpl17}.

\section*{Acknowledgement} JEP acknowledges the support from the UK EPSRC Programme Grant UNLOC EP/J017582/1. JEP and DS are thankful to the Erasmus+ ERC mobility programme between the Aston University and Kharkiv National University that helped us to launch the collaborative activity. AC acknowledges the RISE-FRAMED grant. Discussions with Sergei Turytsyn and Sotos Generalis are gratefully acknowledged.

{\color{OrangeRed}Data used in the figures in this paper are available on \url{http://doi.org/10.17036/researchdata.aston.ac.uk.00000321}.}

\end{document}